\documentclass[a4paper,12pt]{article}
\usepackage[T1]{fontenc}
\usepackage[latin1]{inputenc}
\usepackage{amsthm,amsmath}
\usepackage{amssymb}
\usepackage[all]{xy}
\usepackage{enumerate}
\usepackage{a4wide}
\usepackage{hyperref}
\usepackage{latexsym}
\usepackage{enumerate}
\usepackage[mathscr]{eucal}

\newtheorem{theorem}{Theorem}[section]
\newtheorem{proposition}[theorem]{Proposition}
\newtheorem{corollary}[theorem]{Corollary}
\newtheorem{lemma}[theorem]{Lemma}
\newtheorem*{theorem*}{Theorem}
\newtheorem*{proposition*}{Proposition}
\newtheorem*{corollary*}{Corollary}
\newtheorem*{lemma*}{Lemma}
\theoremstyle{definition}
\newtheorem{definition}[theorem]{Definition}

\newtheorem{example}[theorem]{Example}
\newtheorem{remark}[theorem]{Remark}
\newtheorem*{remark*}{Remark}

\newtheorem*{definition*}{Definition}
\numberwithin{equation}{section}

\newcommand{\w}{{\otimes}}

\newcommand{\T}{{\mathbf{Sub}_{\A,F}(U(b_0))}}

\newcommand{\V}{{\mathcal {V}}}
\newcommand{\G}{{\textbf{ G}}}
\newcommand{\A}{{\mathcal{A}}}
\newcommand{\B}{{\mathcal{B}}}
\newcommand{\C}{{\mathcal{C}}}
\newcommand{\MM}{{\mathcal{M}}}

\newcommand{\e}{{\varepsilon}}
\newcommand{\BB}{{\mathcal{B}_ {\textbf{G}}}}
\newcommand{\K}{{\mathcal {K}}}
\newcommand{\R}{ {\text{Mod}_K}}

\begin{document}
%Title and Authors
\title{A bicategorical version of Masuoka's theorem. Applications to bimodules over functor categories and to firm bimodules\footnote{This research is supported by the grants MTM2004-01406
and MTM2007-61673 from the Ministerio de Edu\-ca\-ci{\'o}n y Ciencia
of Spain, and P06-FQM-01889 from the Consejer{\'\i}a de
Innovaci{\'o}n, Ciencia y Empresa of Andaluc{\'\i}a, with funds from
FEDER (Uni{\'o}n Europea)}}
\author{J. G{\'o}mez-Torrecillas \\
\normalsize Departamento de \'{A}lgebra \\ \normalsize Facultad de Ciencias \\
\normalsize  Universidad de Granada \\ \normalsize E18071 Granada,
Espa\~{n}a
\\ \normalsize
e-mail: \textsf{gomezj@ugr.es} \and B. Mesablishvili \\
\normalsize Mathematical Institute \\ \normalsize Georgian Academy
of Sciences \\  \normalsize
Alexidze Str. 1, \\ \normalsize 0193 Tbilisi, Georgia \\
\normalsize e-mail: \textsf{bachi@rmi.acnet.ge} }

\date{ }
\maketitle
\begin{abstract}We give a bicategorical version of  the main result of   A. Masuoka
({Corings and invertible bimodules,} {\em Tsukuba J. Math.}
\textbf{13} (1989), 353--362) which proposes a non-commutative
version of the fact that for a faithfully flat extension of
commutative rings  $R \subseteq S$,  the relative Picard group
$Pic(S/R)$  is isomorphic to the Amitsur $1$--cohomology group
$H^1(S/R,U)$ with coefficients in the units functor $U$.

\end{abstract}

\section*{Introduction}

In \cite{Masuoka:1989}, A. Masuoka proposed a non-commutative
version of the fact that for a faithfully flat extension of
commutative rings $R \subseteq S$, the relative Picard group
$Pic(S/R)$ is isomorphic to the Amitsur $1$--cohomology group
$H^1(S/R,U)$ with coefficients in the units functor $U$. In the
classical faithfully flat descent setting, each element of
$Pic(S/R)$ is represented by an $R$--subbimodule $J$ of $S$ such
that $S \otimes_R J \cong S$, while the elements of $H^1(S/R,U)$ may
be interpreted as the grouplike elements of the canonical Sweedler's
coring $S \otimes_R S$ or, equivalently, as the endomorphisms of $S
\otimes_R S$ as an $S$--coring. For a general extension of rings $R
\subseteq S$ (unital, but possibly noncommutative), Masuoka's
theorem asserts that, under suitable conditions (e.g., $S_R$
faithfully flat), there exists an isomorphism of monoids $\Gamma :
\mathbf{I}^l_R(S) \rightarrow \mathrm{End}_{S-cor}(S \otimes_R S)$,
where the first one is the monoid of all $R$--subbimodules $J$ of
$S$ such that $S \otimes_R J \cong S$, and the second one consists
of all endomorphisms of the coring $S \otimes_R S$. A more general
version was given in \cite{EG}, where the ring extension $R
\subseteq S$ is replaced by an $R-S$--bimodule $M$ such that $M_S$
is finitely generated and projective, and the role of the Sweedler
canonical coring is now played by the comatrix coring associated to
$M$. In this setting it has been proved in \cite{Me1} that Masuoka's
theorem is valid under the condition of comonadicity of the functor
$ - \otimes_R M$ from right $R$--modules to right $S$--modules.
These generalizations, specially the second one, required of
arguments of categorical nature. Thus, we wonder to what extent a
purely categorical version of Masuoka's theorem was possible. This
paper collects the results of our research on this idea.

In fact, what we prove is a bicategorical version of Masuoka's
theorem, which is rather natural, since both the original result,
and the aforementioned generalizations, can be formulated in the
bicategory of bimodules. Of course, passing from this particular
bicategory to a general one requires to overcome many technical
difficulties, as well as new viewpoints on the objects involved. In
order to make accessible our arguments to readers coming from
different fields, we include rather detailed proofs.

Our generalization of Masuoka's theorem works over an arbitrary
bicategory. Thus, it can be potentially applied to many particular
situations of interest. We have chosen to illustrate how the general
theory works by applying it to two situations close to the case of
bimodules over unital rings, namely, the bicategory of bimodules
over linear categories of functors and the bicategory of firm
bimodules.

The paper is organized as follows. In Section \ref{subcoalgebra} we
consider the comonad generated by a general pair of adjoint functors
and we study the relationship between the possible coalgebra
structures over a given object and some subobjects of its image
under the right adjoint functor. When the left adjoint functor is
comonadic, the aforementioned relationship takes the form of a
bijection, which turns out to be fundamental for our purposes
(Proposition \ref{subobjetoscoalgebrasbiy}).

Section \ref{bicatMasuoka} is the core of the paper. After recalling
some basic notions related to bicategories and monoidal categories,
which serves as well to fix our notations, we proceed to define the
monoid (semigroup with unit)
$\mathbf{I}^l_{\mathcal{V}}(\mathbb{S})$ associated to a monoid
$\mathbb{S}$ in a monoidal category $\mathcal{V}$ (Proposition
\ref{monoide}). When $\mathcal{V}$ is the monoidal category of
bimodules over a ring $R$, the monoid $\mathbb{S}$ is given by a
homomorphism of rings $R \rightarrow S$, and, in the faithfully flat
case, $\mathbf{I}^l_{\mathcal{V}}(\mathbb{S})$ turns out to be
isomorphic to the monoid $\mathbf{I}^l_R(S)$ considered by Masuoka
(see Section \ref{applications}). We make then explicit in
Proposition \ref{reduccionaestricto} that the study of
$\mathbf{I}^l_{\mathcal{V}}(\mathbb{S})$ can be already done in the
strict case. Equipped with these tools, we proceed to the
construction of our generalization of Masuoka's theorem in the
framework of an arbitrary bicategory $\mathbb{B}$. We consider an
adjunction $f \dashv f^* : \mathcal{B} \rightarrow \mathcal{A}$ in
$\mathbb{B}$, the monad (or $\mathcal{A}$--ring) $S_f^{\mathbb{B}}$,
and the comonad (or $\mathcal{B}$--coring)
$\mathfrak{C}_f^{\mathbb{B}}$ built from the adjunction (this is the
general form of a comatrix coring), and we prove (Theorem
\ref{maincomonadic}) that if the functor $\mathbb{B}(\mathcal{A},f):
\mathbb{B}(\mathcal{A},\mathcal{A}) \rightarrow
\mathbb{B}(\mathcal{A},\mathcal{B})$ is comonadic, then it is
possible to define an isomorphism of monoids $\Gamma_f :
\mathbf{I}^l_{\mathbb{B}(\mathcal{A},\mathcal{A})}(S_f^{\mathbb{B}})
\rightarrow
\mathrm{End}_{\mathcal{B}}(\mathfrak{C}_f^{\mathbb{B}},\mathfrak{C}_f^{\mathbb{B}})$,
where the latter stands for the monoid of endomorphisms of the
$\mathcal{B}$--coring $\mathfrak{C}_f^{\mathbb{B}}$. Our strategy to
obtain this isomorphism is to prove it first in the case of
$2$--categories and then deduce the general case. Full details on
this nontrivial process are provided.

In Section \ref{flatpurityseparability} we investigate sufficient
conditions on the $1$--cell $f$ that imply the comonadicity of the
functor $\mathbb{B}(\mathcal{C},f)$, and thus the existence of the
isomorphism of monoids $\Gamma_f$. The nature of these conditions is
flatness and purity (Proposition \ref{flatpure}) or separability
(Proposition \ref{separable}).

In Section \ref{applications} we apply our general results to
bicategories of generalized bimodules. First, we consider the
framework of $K$--linear functor categories, where $K$ is a
commutative ring. These categories are of interest in Representation
Theory of Algebras. We construct a bicategory $\mathbf{Bim}$ of
bimodules whose $0$--cells are small $K$--linear categories. In this
setting, we obtain generalizations of the Masuoka Theorem which boil
down, in the case of $0$--cells of one object (rings) to the results
in \cite{EG} and \cite{Me1}. In particular, we prove that, in the
precomonadic case, the elements of
$\mathbf{I}^l_{\mathbf{Bim}(\mathcal{A},\mathcal{A})}(S_f)$ can be
represented by $\mathcal{A}$--subbimodules of $S_f$, according to
what happens to the case of usual rings and bimodules in \cite{EG}.
The second bicategory to which our general results are applied is
$\mathbf{Firm}$, whose $0$--cells and $1$--cells are, respectively,
firm rings and firm bimodules over them. Firm rings and modules
where considered by D. Quillen in \cite{Quillen:notes} in connection
with the excision problem in homology. We prove that $\mathbf{Firm}$
is a biclosed bicategory, which allow the characterization of those
$1$--cells that have a right adjoint. We get then a version of
Masuoka's theorem for firm bimodules which admit right adjoint in
$\mathbf{Firm}$.

\section{Subobjects and coalgebra structures}\label{subcoalgebra}

Let $\A$ be a category and let $ a_0 \in \A$. We write
$\text{Sub}_{\A}(a_0)$ for the class of all monomorphisms with
codomain $a_0$. Recall that monomorphisms $i_a: a \to a_0$ in $\A$
are pre-ordered by setting $i_a \leq i_{a'}$ if $i_a=i_{a'} k$ for
some $k: a \to a'$; two such monomorphisms are equivalent if $i_a
\leq i_{a'}$ and $i_{a'} \leq i_a$ (i.e., if there exists an
isomorphism $f: a \to a'$ with $i_a= i_{a'}f$), and the equivalence
classes are called subobjects of $a_0$. If $i_a: a \to a_0$ is a
monomorphism in $\A$, we write $[i_a]$ for the corresponding
subobject of $a_0$. We will say that an object $x$ of $\A$ is a
subobject of $a_0$ if there is a monomorphism $x \to a_0$. We let
$\text{\textbf{Sub}}_{\A}(a_0)$ denote the class of subobjects of
$a_0$.

If $ \eta, \epsilon : F \dashv U : \B \to \A$ is an adjunction with
adjunction isomorphism $$\alpha_{x,y}: \A (x, U(y)) \to \B(F(x),
y)$$ and if $b_0 \in \B$, write $\text{\textbf{Sub}}_{\A,F}(U(b_0))$
for the subclass of $\text{\textbf{Sub}}_{\A}(U(b_0))$ whose
elements  are those subobjects $ [i_a : a \to U(b_0)]$ for which the
morphism
$$\xi^F_{i_a}=\alpha_{a, \,b_0}(i_a)=\varepsilon_{b_0} \cdot F(i_a): F(a) \to b_0$$
is an isomorphism. Note that $\text{\textbf{Sub}}_{\A,\,F}(U(b_0))$
is well- defined: if $ [i_{a}]=[i_{a'}]$ in
$\text{\textbf{Sub}}_{\A}(U(b_0))$, then there exists an isomorphism
$f: a \to a'$ such that $i_{a} \cdot f=i_{a'}$, and one has
$$\xi^F_{i_{a'}}=\alpha_{a'\!, \,b_0} (i_{a'})=\varepsilon_{b_0}
\cdot F(i_{a'})=\varepsilon_{b_0} \cdot F(i_a) \cdot F(f)=\alpha_{a,
b_0} (i_a) \cdot F(f)=\xi^F_{i_a} \cdot F(f),$$ and hence
$\xi^F_{i_{a'}}=\alpha_{a'\!,\, b_0} (i_{a'})$ is an isomorphism
too. We will generally drop the $\A$ (resp. $F$) from the notation
$\text{\textbf{Sub}}_{\A,\,F}(U(b_0))$ (resp. $\xi^F_{i_a}$) when
there is no danger of confusion.

We now fix an adjunction $ \eta, \epsilon : F \dashv U : \B \to \A$
with adjunction isomorphism $$\alpha_{x,\, y}: \A (x, U(y)) \to
\B(F(x), y)$$ and an object $b_0 \in \B$.

\begin{lemma}\label{conservativosubobjetos}
Suppose that $F$ is conservative (=isomorphism reflecting). If
\,$[i_a]$ and $[i_{a'}]$ are elements of
$\mathbf{Sub}_{\A,\,F}(U(b_0))$ such that $i_a \leq i_{a'}$ in
$\mathbf{Sub}_{\A}(U(b_0))$, then $[i_a]=[i_{a'}].$
\end{lemma}
\begin{proof} Since $i_a \leq i_{a'}$ in $\mathbf{Sub}_{\A}(U(b_0))$,
there exists a morphism $f: a \to a'$ such that $i_a=i_{a'} \cdot
f$, and then we have $$\xi^F_{i_{a'}}\cdot
F(f)=\alpha_{a',\,b_0}\cdot(i_{a'})\cdot F(f)=\e_{b_0} \cdot
F(i_{a'})\cdot F(f)=\e_{b_0} \cdot F(i_{a'}f)=\e_{b_0} \cdot U(i_a)=
\alpha_{a,\,b_0}(i_a)=\xi^F_{i_a}.$$ It follows that $F(f)$ is an
isomorphism and thus $f$ is also an isomorphism since $F$ is assumed
to be conservative. Thus $[i_a]=[i_{a'}].$
\end{proof}

In order to proceed, we need the following
\begin{lemma}\label{lema2}
A morphism $f: F(a) \to F(a')$ is an isomorphism if and only if the
morphism $U(f): UF(a) \to UF(a')$ is so.
\end{lemma}
\begin{proof}One direction is clear, so suppose that $U(f): UF(a) \to UF(a')$
is an isomorphism. In the following diagram
$$
\xymatrix { FUFUF(a)
\ar[dd]_{FUFU(f)}\ar@{->}@<0.5ex>[rr]^-{FU\varepsilon_{F(a)}} \ar@
{->}@<-0.5ex> [rr]_-{\varepsilon_{FUF(a)}} && FUF(a) \ar[dd]^{FU(f)}
\ar@{->}@<0ex>[rr]^-{\varepsilon_{F(a)}} && F(a)\ar[dd]^{f}\\\\
FUFUF(a') \ar@{->}@<0.5ex>[rr]^-{FU\varepsilon_{F(a')}} \ar@
{->}@<-0.5ex> [rr]_-{\varepsilon_{FUF(a')}} && FUF(a')
\ar@{->}@<0ex>[rr]_-{\varepsilon_{F(a')}} && F(a') ,}
$$ which is commutative by naturality of $\varepsilon$, each row
is a (split) coequalizer (see \cite{M}), implying -since $FUFU(f)$
and $FU(f)$ are isomorphisms- that $f$ is an isomorphism.
\end{proof}

Let us write $\mathbf{T}$ for the monad on $\A$ generated by the
adjunction $F \dashv U : \B \to \A$, $\eta^{\textbf{T}},
\varepsilon^{\textbf{T}}: F^{\textbf{T}} \dashv U^{\textbf{T}} :
\A^{\textbf{T}} \to \A$ for the corresponding free-forgetful
adjunction and consider an arbitrary object $a_0$ of $\A$. Since
$UF=U^{\textbf{T}}F^{\textbf{T}}$, clearly
$\text{\textbf{Sub}}_{\A}(UF(a_0))=\text{\textbf{Sub}}_{\A}(U^{\textbf{T}}F^{\textbf{T}}(a_0))$.

\begin{proposition}\label{subobjetosiguales}
 For any $[(a, i_a)]\in
\mathbf{Sub}_{\A}(UF(a_0))$, $[(a, i_a)]\in
\mathbf{Sub}_{\A,\,F}(UF(a_0))$ if and only if $[(a, i_a)]\in
\mathbf{Sub}_{\A,\,F^{\mathbf{T}}}(U^{\mathbf{T}}(a_0))$.
\end{proposition}
\begin{proof} Writing $K^{\textbf{T}}: \B \to \A^{\textbf{T}}$ for
the comparison functor, one sees -since $K^{\textbf{T}}F
=F^{\textbf{T}}$ and $\varepsilon^{\textbf{T}}
K^{\textbf{T}}=K^{\textbf{T}}\varepsilon$ (see, \cite{M})- that the
morphism
$$
\xymatrix{ \xi^{F^{\textbf{T}}}_{i_a}: F^{\textbf{T}}(a)
\ar[rr]^-{F^{\textbf{T}}(i_a)}&&
F^{\textbf{T}}U^{\textbf{T}}F^{\textbf{T}}(a_0)
\ar[rr]^-{\varepsilon^{\textbf{T}}_{F^{\textbf{T}}(a_0)}} &&
F^{\textbf{T}}(a_0)},
$$ can be rewritten as the composite
$$
\xymatrix{K^{\textbf{T}} F(a) \ar[rr]^-{K^{\textbf{T}}F(i_a)}&&
K^{\textbf{T}}FUF(a_0) \ar[rr]^-{K^{\textbf{T}}\varepsilon_{F(a_0)}}
&& K^{\textbf{T}}F(a_0)}.
$$ It is now easy to see that
$$U^{\textbf{T}}(\xi^{F^{\textbf{T}}}_{i_a})=U\varepsilon_{F(a_0)}\cdot
UF(i_a)=U(\varepsilon_{F(a_0)}\cdot F(i_a))=U(\xi^F_{i_a}).$$
Applying now Lemma \ref{lema2} we get that $\xi^F_{i_a}$ is an
isomorphism if and only if
$U^{\textbf{T}}(\xi^{F^{\textbf{T}}}_{i_a})$ (and hence
$\xi^{F^{\textbf{T}}}_{i_a}$, since $U^{\textbf{T}}$ preserves and
reflects isomorphisms) is so.
\end{proof}

As an easy consequence, we have

\begin{proposition}\label{submon}
In the situation of Proposition \ref{subobjetosiguales},
$$\mathbf{Sub}_{\A,\,F}(UF(a_0))=
\mathbf{Sub}_{\A,\,F^{\mathbf{T}}}(U^{\mathbf{T}}F^{\mathbf{T}}(a_0)).$$
\end{proposition}

Given an arbitrary monad $\textbf{T}=(T, \eta, \mu)$ on $\A$ and an
object $a_0 \in \A$, we write $\text{\textbf{Sub}}_{\A,
\textbf{T}}(T(a_0))$ for the subclass of
$\text{\textbf{Sub}}_{\A}(T(a_0))$ whose elements are those
subobjects $[(a, i_a)]$ for which the composite
$$\xymatrix{
\xi^{\textbf{T}}_{i_a}: T(a) \ar[r]^{T(i_a)}& T^2 (a_0)
\ar[r]^{\mu_{a_0}}& T(a_0)}
$$ is an isomorphism.

\begin{proposition}\label{mon1}
\begin{itemize}
\item[\emph{(i)}]Let $\mathbf{T}$  be a monad on a category $\A$.
If $f:a_0 \to a_0'$ is an isomorphism in $\A$, then the assignment
$$(i_a: a \to a_0)\longrightarrow (fi_a:a \to a_0')$$ yields a
bijection $$\mathbf{Sub}_{\A,
\mathbf{T}}(T(a_0))\simeq\mathbf{Sub}_{\A, \mathbf{T}}(T(a'_0)).$$

\item[\emph{(ii)}]Let $\mathbf{T}$ and $\mathbf{T}'$ be monads on
a category $\A$. If $\tau: \mathbf{T} \to \mathbf{T}' $ is an
isomorphism of monads, then for any $a_0 \in \A$, the assignment
$$[(a, i_a: a \to T(a_0))]\longrightarrow [(a, \tau_{a_0}\cdot i_a: a \to T'(a_0))]
$$ yields a bijection $$\mathbf{Sub}_{\A,
\mathbf{T}}(T(a_0))\simeq\mathbf{Sub}_{\A, \mathbf{T}'}(T'(a_0)).$$
\end{itemize}
\end{proposition}

\begin{proposition}For any monad $\mathbf{T}=(T, \eta, \mu)$ on $\A$ and
an object $a_0 \in \A$,
$$\mathbf{Sub}_{\A, \mathbf{T}}(T(a_0))=\mathbf{Sub}
_{\A,\,F^{\mathbf{T}}}(U^{\mathbf{T}}(a_0)),$$ where $F^{\mathbf{T}}
: \A \to \A^{\mathbf{T}}$ is the free $\mathbf{T}$-algebra functor.
\end{proposition}
\begin{proof} Since the forgetful functor $U^{\textbf{T}} : \A^{\textbf{T}} \to
\A$ preserves and reflects isomorphisms, the morphism
$\xi^{F^{\textbf{T}}}_{i_a}$ is an isomorphism if and only if
$U^{\textbf{T}}(\xi^{F^{\textbf{T}}}_{i_a})$ is so. But since
$U^{\textbf{T}} F^{\textbf{T}}=T,$ it is easy to check that
$U^{\textbf{T}}(\xi^{F^{\textbf{T}}}_{i_a})=\xi^{\textbf{T}}_{i_a}$,
proving that $\xi^{F^{\textbf{T}}}_{i_a}$ is an isomorphism if and
only if $\xi^{\textbf{T}}_{i_a}$ is so.
\end{proof}

Combining this with Proposition \ref{submon}, we get

\begin{proposition}\label{mon2}
Let $\mathbf{T}=(UF, \eta, U\varepsilon F)$ be the monad
corresponding to the adjunction $\eta, \varepsilon : F \dashv U : \B
\to \A$, then for any object $a_0 \in \A$,
$$\mathbf{Sub}_{\A, \mathbf{T}}(T(a_0))=\mathbf{Sub}_{\A,\,F}(UF(a_0))=
\mathbf{Sub}_{\A,\,F^{\mathbf{T}}}(U^{\mathbf{T}}F^{\mathbf{T}}(a_0)).$$
\end{proposition}

\bigskip

Let $\G=\!\!(G=\!FU, {\e_{\G}}=\e, \delta_{\G}=F\eta U)$ be the
comonad on $\B$ generated by the adjunction $F \dashv U$,
$\B_{\!\textbf{G}}$ the corresponding Eilenberg-Moore category of
$\G$-coalgebras, $U_{\!\G}: \B_{\!\textbf{G}} \to \B$ the forgetful
functor, and $K_{\!\G}: \A \to \B_{\!\G}$ the comparison functor
that takes an arbitrary object $a \in \A$ to the
$\textbf{G}$-coalgebra $(F(a), F(\eta_a))$. We write
$\textbf{G-}\text{coalg}(b_0)$ for the set of all $\G$-coalgebra
structures on $b_0$.

\begin{proposition}The assignment $$\xymatrix {
[(a,i_a)] \ar[r]& (b_0 \ar[r]^{\xi^{-1}_a}& F(a) \ar[r]^-{F(\eta_a
)}& FUF(a) \ar[r]^-{FU(\xi_a)}& FU(b_0)=G(b_0))}$$ defines a map
$$\Psi_{F,b_0} :\mathbf{Sub}_{\A,F}(U(b_0)) \to \mathbf{G-}\text{coalg}(b_0).$$
\end{proposition}

\begin{proof}Let $[(a,i_a)] \in \text{\textbf{Sub}}_{\A,F}(U(b_0))$. We have to
show that ($\Psi_{F,b_0}((a,\xi_a))$ is a $\G$-coalgebra morphism
-to be replaced by) $\Psi_{F,b_0}(a,\xi_a)\in
\mathbf{G-}\text{coalg}(b_0)$ . First we prove that
$\varepsilon_{b_0} \cdot \Psi_{F,b_0}((a,\xi_a))=1.$ Using the
naturality of $\varepsilon : FU \to 1$ and the triangular identity
$\varepsilon_{F(a)} \cdot F(\eta_a)=1_{F(a)}$ for the adjunction
$\eta, \epsilon : F\dashv U$, we see that $$\varepsilon _{b_0} \cdot
\Psi_{F,b_0}((a,\xi_a))=\varepsilon _{b_0} \cdot FU(\xi_a) \cdot
F(\eta_a) \cdot \xi^{-1}_a=$$$$=\xi_a \cdot \varepsilon_{F(a)} \cdot
F(\eta_a)\cdot \xi^{-1}_a=\xi_a \cdot 1_{F(a)}\cdot \xi^{-1}_a=\xi_a
\cdot \xi^{-1}_a=1_{b_0}.$$ Next, we have to prove that
$$(\delta_{\textbf{G}}){_{b_0}} \cdot \Psi_{F,b_0}
((a,\xi_a))= G(\Psi_{F,b_0}((a,\xi_a))) \cdot
\Psi_{F,b_0}((a,\xi_a)).$$ Since $\eta : 1 \to UF$ is a natural
transformation, the diagrams $$\xymatrix{ a \ar[r]^{\eta_a}
\ar[d]_{\eta_a} & UF(a) \ar[d]^{UF(\eta_a)} & \text{and}&
UF(a)\ar[d]_{U(\xi_a)}
\ar[r]^-{\eta_{UF(a)}} & UFUF(a) \ar[d]^{UFU(\xi_a)}\\
UF(a) \ar[r]_-{\eta_{UF(a)}} & UFUF(a) && U(b_0)
\ar[r]_-{\eta_{U(b_0)}} & UFU(b) }
$$ commute. Hence

\begin{equation}
FUF(\eta_a) \cdot F(\eta_a)=F(\eta_{UF(a)})\cdot F(\eta_a)
\end{equation}
and
\begin{equation}
FUFU(\xi_a)\cdot F(\eta_{UF(a)})=F(\eta_{U(b_0)})\cdot FU(\xi_a).
\end{equation} We now have:

$$G(\Psi_{F,b_0}((a,\xi_a)))
\cdot \Psi_{F,b_0}((a,\xi_a))=$$$$=FUFU(\xi_a)\cdot FUF(\eta_a)
\cdot FU(\xi^{-1}_a) \cdot FU(\xi_a) \cdot F(\eta_a) \cdot
\xi^{-1}_a=$$$$=FUFU(\xi_a)\cdot FUF(\eta_a) \cdot F(\eta_a) \cdot
\xi^{-1}_a= \,\,\,\,\,\, \text{by (1)}$$$$=FUFU(\xi_a)\cdot
F(\eta_{UF(a)})\cdot F(\eta_a) \cdot \xi^{-1}_a= \,\,\,\,\,\,
\text{by (2)}$$$$=F(\eta_{U(b_0)})\cdot FU(\xi_a)\cdot F(\eta_a)
\cdot \xi^{-1}_a=$$$$=F(\eta_{U(b_0)})\cdot
\Psi_{F,b_0}((a,\xi_a))=$$$$(\text{since}\,\,\,\,\,\,
(\delta_{\textbf{G}}){_{b_0}}=F(\eta_{U(b_0)}) \,\,\,\,\,\, \text{by
definition of}\,\,\,\,\,\,  \delta_{\textbf{G}})$$
$$=(\delta_{\textbf{G}}){_{b_0}}\cdot \Psi_{F,b_0}((a,\xi_a)).$$

Suppose now that $[\!(a, \xi_a)\!]=[(a',\xi_{a'})] $ in
$\text{\textbf{Sub}}_{\A,F}(U(b_0))$. Then there exists an
isomorphism $f:a \to a'$ such that $i_{a'}f=i_a$, and thus
$\xi_{a'}F(f)=\xi_a$. Then, since $\eta : 1 \to UF$ is a natural
transformation, the diagram
$$\xymatrix {
b_0 \ar[rdd]_{\xi^{-1}_{a'}}\ar[r]^{\xi^{-1}_a} &F(a)
\ar[rr]^{F(\eta_a)} && FUF(a)
\ar[rr]^{FU(\xi_a)}&& FU(b_0)\\\\
& F(a') \ar[uu]_{F(f^{-1})} \ar[rr]_{F(\eta_{a'})}&&
FUF(a_1)\ar[uu]_{FUF(f^{-1})} \ar[rr]_{FU(\xi_{a'})}&& FUF(b_0)
\ar@{=}[uu]}
$$ commutes. It follows that
$${\Psi}_{F,b_0}([(a,i_a)])=FU(\xi_a) \cdot F(\eta_a) \cdot
\xi^{-1}_a=$$$$=FU(\xi_{a'}) \cdot F(\eta_{a'}) \cdot
\xi^{-1}_{a'}={\Psi}_{F,b_0}([(a',i_{a'})]).$$ Thus, $\Psi_{F,b_0}$
is well-defined.
\end{proof}

Recall that $F$ is called (pre)comonadic if the comparison functor
$K_{\textbf{G}}: \A \to \B_{\textbf{G}}$ is an equivalence of
categories (full and faithful). We are ready to state the main
result of this section.

\begin{proposition}\label{subobjetoscoalgebrasbiy}
If $F$ is comonadic, the map $$\Psi_{F,b_0}:
 \T \to \mathbf{G-}\text{coalg}(b_0)$$ is bijective.
\end{proposition}

\begin{proof} We first show that $\Psi_{F,b_0}$ is injective. Indeed,
suppose that $\Psi_{F,b_0}([\!(a,i_a\!])=
\Psi_{F,b_0}([\!(a',i_{a'})\!]).$ Then the diagram
$$\xymatrix {
& F(a) \ar[rr]^{F(\eta_a)}&& FUF(a) \ar[rd]^{FU(\xi_a)} &\\
b_0 \ar[dr]_{\xi^{-1}_{a'}}\ar[ru]^{\xi^{-1}_a} &&&& FU(b_0)\\
& F(a') \ar[rr]_{F(\eta_{a'})}&& FUF(a') \ar[ru]_{FU(\xi_{a'})} &}
$$ commutes. It follows that we have the following commutative diagram
$$\xymatrix {
F(a) \ar[dd]_{\xi^{-1}_{a'} \cdot \xi_{a}}\ar[rr]^{F(\eta_a)} &&
FUF(a) \ar[dd]^{FU(\xi^{-1}_{a'} \cdot
\xi_{a})}\\\\
F(a') \ar[rr]_{F(\eta_{a'})}&& FUF(a').}
$$ But to say that this diagram commutes is just to say that
$\xi^{-1}_{a'} \cdot \xi_{a}:F(a)\to F(a')$ is a morphism in $\BB$
from the $\G$-coalgebra $K_{\!\G}(a)=(F(a), F(\eta_a))$ to the
$\G$-coalgebra $K_{\!\G}(a')=(F(a'), F(\eta_{a'}))$ and since the
functor $F$ is assumed to be comonadic, the functor $K_{\!\G}$ is an
equivalence of categories and thus there exits a morphism $f: a \to
a'$ in $\A$ such that $F(f)=\xi^{-1}_{a'} \cdot \xi_{a}$, whence
$\xi_a = \xi_{a'} \cdot F(f)$. It follows that $F(f)$ is an
isomorphism and hence so also is $f$, since $F$ is comonadic and any
comonadic functor reflects isomorphisms. Moreover, since
$\alpha_{a,\,b_0}(i_a)=\xi_{i_a}= \xi_{i_{a'}}\cdot
F(f)=\alpha_{a',\,b_0}(i_{a'})\cdot F(f)=\e_{b_0} \cdot F(i_{a'})
\cdot F(f)=\e_{b_0}\cdot F(i_{a'}\cdot f)= \alpha_{a,\,b_0}
(i_{a'}\cdot f)$ and since $\alpha_{a,\, b_0}$ is bijective, one
concludes that $i_{a'}f=i_a.$ Thus, $[\!(a,i_a)\!]=[\!(a',
i_{a'})\!]$, i.e., $\Psi_{F,b_0}$ is injective.

Suppose now that $\theta_{b_0} : b_0 \to G(b_0)$ is such that
$(b_0,\theta_{b_0}) \!\in \! \mathbf{G-}\text{coalg}(b_0)$, and
consider $\tilde{b_0}=R_{\! \G}(b_0,\theta_{b_0})$, where $R_{\! \G}
: \B_{\!\G} \to \A$ is the right adjoint of the comparison functor
$K_{\!\G} : \A \to \B_{\!\G}$. It is well known (see \cite{M}) that
$\tilde{b}_0$ appears as the equalizer
$$\xymatrix {\tilde{b}_0 \ar[r]^-{e_{\tilde{b}_0}} & U(b_0)
\ar@{->}@<0.5ex>[r]^-{U(\theta_{b_0})} \ar@ {->}@<-0.5ex>
[r]_-{\eta_{U(b_0)}}& UFU(b_0).}$$ Hence
$(\tilde{b_0},e_{\tilde{b}_0}) \in \text{Sub}_{\A}(U(b_0)).$ If
$\varepsilon_1 :K_{\!\G}  R_{\!\G} \to 1$ is the counit of the
adjunction $K_{\!\G} \dashv R_{\!\G}$ (which is an isomorphism,
since $F$ is assumed to be comonadic), then the diagram
$$\xymatrix {
F(\tilde{b}_0) \ar[d]_{(\varepsilon_1)_{(b_0, \theta_{b_0})}}\ar[r]^
{F(e_{\tilde{b}_0})}& FU(b_0)\\
b_0 \ar[ru]_{\theta_{b_0}},}
$$ where $(\varepsilon_1)_{(b_0, \theta_{b_0})}$ denotes the $(b_0,
\theta_{b_0})$-component of $\varepsilon_1$, commutes (see, for
example, \cite{B}). It follows that
\begin{equation}
\theta_{b_0}=F(e_{\tilde{b}_0}) \cdot ((\varepsilon_1)_{(b_0,\,
\theta_{b_0})})^{-1}
\end{equation}
and that
\begin{equation}
\varepsilon_{b_0} \cdot F(e_{\tilde{b}_0})= \varepsilon_{b_0} \cdot
\theta_{b_0} \cdot (\varepsilon_1)_{(b_0,
\theta_{b_0})}=(\varepsilon_1)_{(b_0, \theta_{b_0})}.
\end{equation} Hence $\xi_{e_{\tilde{b}_0}}=\alpha_{\tilde{b}_0, b_0}
(e_{\tilde{b}_0})=\varepsilon_{b_0} \cdot F(e_{\tilde{b}_0})=
(\varepsilon_1)_{(b_0, \theta_{b_0})}$ is an isomorphism. It follows
that the assignment $$(b_0, \theta_{b_0}) \to [(\tilde{b}_0,
e_{\tilde{b}_0})]$$ defines a map $$\overline{\Psi}_{F, b_0}:
\textbf{G-}\text{coalg}(b_0)) \to
\text{\textbf{Sub}}_{\A,F}(U(b_0)).$$

We now claim that $\Psi_{F,b_0} \cdot \bar{\Psi}_{F,b_0}=1$. Indeed,
we have
$$\Psi_{F,b_0}(\bar{\Psi}_{F,b_0}(b_0,
\theta_{b_0})))=\Psi_{F,b_0}(\tilde{b}_0,e_{\tilde{b}_0})=
FU(\xi_{e_{\tilde{b}_0}})\cdot F(\eta_{\tilde{b}_0})\cdot
(\xi_{e_{\tilde{b}_0}})^{-1} =$$$$=FU((\varepsilon_1)_{(b_0,
\theta_{b_0}}) \cdot F(\eta_{\tilde{b}_0}) \cdot
((\varepsilon_1)_{(b_0, \theta_{b_0})})^{-1}=\,\,\,\,\, \text{by
(4)}
$$$$=FU(\varepsilon_{b_0})\cdot FUF(e_{\tilde{b}_0}) \cdot F
(\eta_{\tilde{b}_0}) \cdot ((\varepsilon_1)_{(b_0,
\theta_{b_0})})^{-1}=$$$$\,\,\,\,\, \text{(by the naturality of}\,\,
\eta)$$$$=FU(\varepsilon_{b_0})\cdot F(\eta_{U(b_0)}) \cdot
F(e_{\tilde{b}_0}) \cdot ((\varepsilon_1)_{(b_0,
\theta_{b_0})})^{-1}=\,\,\,\,\,\,\,\,\,\,\,\,\,\,\,\,\,\,$$
$$=F(U(\varepsilon_{b_0}))\cdot \eta_{U(b_0))} \cdot
F(e_{\tilde{b}_0}) \cdot ((\varepsilon_1)_{(b_0,
\theta_{b_0})})^{-1}=\,\,\,\,\,\,\,\,\,\,\,\,$$$$\,\,\, \text{(by
the triangular identity)}$$$$=F(e_{\tilde{b}_0}) \cdot
((\varepsilon_1)_{(b_0, \theta_{b_0})})^{-1}=\,\,\,\,\,\, \text{by
(3)}$$$$=\theta _{b_0}.$$ So that $\Psi_{F,b_0} \cdot
\bar{\Psi}_{F,b_0}=1$, and since $\Psi_{F,b_0}$ is injective, it is
bijective. This completes the proof.
\end{proof}

We close this section showing that Proposition
\ref{subobjetoscoalgebrasbiy} might be viewed as a categorification
of \cite[Theorem 1.2]{Nuss/Wambst:2007}.

\begin{example}
Let $R \subseteq S$ be an extension of unital rings, and consider
the associated adjunction $\eta, \epsilon: - \otimes_R S \dashv U :
\mathrm{Mod}_S \rightarrow \mathrm{Mod}_R$, where $U$ denotes the
restriction of scalars functor from right $S$--modules to right
$R$--modules. The corresponding comonad $\mathbf{G}$ is isomorphic
to the one given by the Sweedler canonical $S$--coring $S \otimes_R
S$, that is, $\mathbf{G} \cong - \otimes_S S \otimes_R S :
\mathrm{Mod}_S \rightarrow \mathrm{Mod}_S$. Thus coalgebras over
$\mathbf{G}$ are identified with right $S \otimes_R S$--comodules.
Therefore, for a given right $S$--module $M_S$, $\mathbf{G}-coalg
(M)$ is the set of all right $S \otimes_R S$--coactions (i.e.,
comodule structures) on $M_S$. Therefore, $\mathbf{G}-coalg (M) =
Z^1(S \otimes_R S, M)$, the set of descent $1$--cocycles on $S
\otimes_R S$ with values in $M$, according to
\cite{Brzezinski:2008}. Now, let us interpret the set
$\mathbf{Sub}_{\mathrm{Mod}_R,- \otimes_R S}(U(M))$. From the
definition, this is the set of all right $R$--submodules $N$ of
$M_R$ such that the map
\[
N \otimes_R S \rightarrow M \quad ( n \otimes_R s \mapsto ns)
\]
is bijective. If $M = N_0 \otimes_R S$ for a given right $R$--module
$N_0$, then each $N \in \mathbf{Sub}_{\mathrm{Mod}_R,- \otimes_R
S}(U(M))$ becomes a \emph{twisted form of $N_0$ over $S/R$}
according to \cite{Nuss/Wambst:2007}. In fact, if we denote the set
of such twisted forms by $\mathrm{twist}(S/R,N_0)$, then
\[
\mathbf{Sub}_{\mathrm{Mod}_R,- \otimes_R S}(U(N_0 \otimes_R S)) = \{
N \in \mathrm{twist}(S/R,N_0) \,| \, \eta_N \text{ is monic } \}
\]
Hence, when $\eta_N$ is a monomorphism for all $N \in
\mathrm{Mod}_R$ (e.g., if $- \otimes_R S : \mathrm{Mod}_R
\rightarrow \mathrm{Mod}_R$ is comonadic), we get
\[
\mathbf{Sub}_{\mathrm{Mod}_R,- \otimes_R S}(U(N_0 \otimes_R S)) =
\mathrm{twist}(S/R,N_0),
\]
and, therefore, Proposition \ref{subobjetoscoalgebrasbiy} gives in
particular a bijection
\begin{equation}\label{twistcocycles}
\mathrm{twist}(S/R,N_0) \cong Z^1(S \otimes_R S, N_0 \otimes_R S).
\end{equation}
When $R \subseteq S$ is an $H$-Hopf-Galois extension, where $H$ is a
Hopf algebra, then, by using the well-known interplay between Hopf
$S$--modules and comodules over the $S$--coring $S \otimes_R S$
(see, e.g., \cite[2.6]{Brzezinski:2008}), one easily deduces from
\eqref{twistcocycles} the bijection
\[
\mathrm{Twist}(S/R,N_0) \cong  H^1(H,N_0 \otimes_R S)
\]
proved in \cite[Theorem 1.2]{Nuss/Wambst:2007}, where $H^1(H,M)$ is
the first (nonabelian) cohomology set of the Hopf algebra $H$ with
coefficients in the $(H,S)$--Hopf module $N_0 \otimes S$ (as defined
in \cite{Nuss/Wambst:2007}).
\end{example}

\bigskip

\section{A bicategorical approach to Masuoka's
theorem}\label{bicatMasuoka}

Our aim in this section is to formulate and prove a version of
Masuoka's theorem for a general bicategory. We begin by recalling
that a bicategory $\mathbb{B}$ consists of :
\begin{itemize}
\item a class $\text{Ob}(\mathbb{B})$ of objects, or 0-cells;

\item a family $\mathbb{B}(\A, \B)$, for all $\A, \B \in
\text{Ob}(\mathbb{B})$, of hom-categories, whose objects and
morphisms are respectively called 1-cells and 2-cells;

\item a (horizontal) composition operation, given by a family of
functors
$$\mathbb{B}(\B, \C) \times \mathbb{B}(\A, \B) \to \mathbb{B}(\A,
\C)$$ whose action on a pair $(g, f) \in \mathbb{B}(\B, \C)\times
\mathbb{B}(\A, \B)$ is written $g \otimes f$;

\item identities, given by 1-cells $1_\A \in \mathbb{B}(\A, \A)$,
for $\A \in \text{Ob}(\mathbb{B})$;

\item natural isomorphisms $$\alpha_{h, g, f}: (h \otimes g)
\otimes f \simeq h \otimes (g \otimes f), l_f : 1_\A \otimes f
\simeq f \,\,\text{and} \,\, r_f : f \otimes 1_\A \simeq f, $$
\end{itemize} subject to three coherence axioms (see \cite{Bn}).

\medskip

As a basic example, consider the bicategory $\mathsf{Bim}$ whose
$0$--cells are associative rings with unit, and, for each pair of
rings $A, B$, the hom-category $\mathsf{Bim}(A,B)$ is the category
of all unital $A$--$B$--bimodules and all homomomorphisms of
$A$--$B$--bimodules between them. The horizontal composition is the
opposite of the usual tensor product of bimodules.

\medskip

A \emph{strict} bicategory, or 2-category is a bicategory in which
$\alpha, l$ and $r$ are all identities. When $\mathbb{B}$ is a
2-category, then the composition operation sign $\otimes$ is usually
omitted.

\medskip

An example of $2$--category is $\mathrm{CAT}$, with categories as
$0$--cells, functors as $1$--cells and natural transformations as
$2$--cells.

\medskip

For each bicategory $\mathbb{B}$, there exists the \emph{transpose}
bicategory $\mathbb{B}^t$, defined by:
\begin{itemize}
\item $\text{Ob}(\mathbb{B}^t)=\text{Ob}(\mathbb{B})$,

\item $\mathbb{B}^t(\A, \B)=\mathbb{B}(\B, \A)$ for all objects
$\A, \B \in \text{Ob}(\mathbb{B})$,

\item the horizontal composition operation $$-\otimes^t - :
\mathbb{B}^t(\B, \C) \times \mathbb{B}^t(\A, \B)\to \mathbb{B}^t(\A,
\C)$$ is given by the composite
$$
\xymatrix{\mathbb{B}^t(\B, \C) \times \mathbb{B}^t(\A,
\B)=\mathbb{B}(\C, \B) \times \mathbb{B}(\B, \A)\simeq
\mathbb{B}(\B, \A) \times \mathbb{B}(\C, \B) \ar[r]^-{-\otimes -}&
\mathbb{B}(\C, \A)=\mathbb{B}^t(\A, \C),}$$

\item $1^t_\A=1_\A$,

\item $l^t=r$,

\item $r^t =l$,

\item $\alpha^t =\alpha^{-1}$.
\end{itemize}

A \emph{homomorphism } $\Phi=(\phi, \overline{\phi}, \phi_0):
\mathbb{B} \to \mathbb{B'} $ between bicategories consists of:

\begin{itemize}
\item a function $\phi: \text{Ob}(\mathbb{B}) \to
\text{Ob}(\mathbb{B'})$;

\item functors $\phi_{\A, \B}: \mathbb{B}(\A, \B) \to
\mathbb{B'}(\phi(\A), \phi(\B))$;

\item natural isomorphisms $$\overline{\phi}_{f, g}: \phi(g)
\otimes' \phi(f) \to \phi (g \otimes f)$$ and $$(\phi_0)_\A :
\phi(1_\A) \to \phi(1_\A),$$
\end{itemize} subject to appropriate coherence axioms (see \cite{Bn}).

Given two homomorphisms $\Phi, \Phi' : \mathbb{B} \to \mathbb{B'}$
with $\phi(\A)=\phi'(\A)$ for all objects $\A \in
\text{Ob}(\mathbb{B})$, an \emph{icon} $\tau: \Phi \to \Phi'$ from
$\Phi$ to $\Phi'$ consists of a natural transformation
$$\tau_{\A, \B}: \phi_{\A, \B} \to \phi'_{\A, \B} :\mathbb{B}(\A, \B)\to
\mathbb{B'}(\phi(\A), \phi(\B))=\mathbb{B'}(\phi'(\A), \phi'(\B))$$
for all objects $\A, \B \in \text{Ob}(\mathbb{B})$ satisfying two
coherence conditions- one for horizontal composition, one for
identity 1-cells (see \cite{L}).

\medskip

Recall that a monoidal category is a one-object bicategory, that is,
a monoidal category $\V$ is an ordinary category $V$ with a
bifunctor
$$-\otimes - : V \times V \to V$$ and a unit object $I \in V$,
together with natural isomorphisms
$$\alpha_{X,Y, Z}: (X \otimes Y )\otimes Z \simeq X \otimes (Y \otimes Z),$$
$$l_X: I \otimes X \simeq X,$$ and $$r_X: X \otimes I \simeq X,$$
subject to appropriate coherence axioms (see, for example,
\cite{M}). Recall also that a monoidal category $\V=(V, \otimes, I,
\alpha, l, r)$ is said to be \emph{strict} if the structure
morphisms $\alpha, l, r$ are all identities.

Note that if $\mathbb{B}$ is a bicategory, then for any $\A \in
\text{Ob}(\mathbb{B})$, $(\mathbb{B}(\A, \A), \otimes, 1_\A, \alpha,
l, r)$ is a monoidal category.

\medskip

 A monoidal functor $$\Phi=(\phi, \overline{\phi}, \phi_0):
\V=(V, \otimes, I, \alpha, l, r) \to \V'=(V', \otimes', I', \alpha',
l', r')$$ between monoidal categories consists of:
\begin{itemize}
\item an ordinary functor $\phi: V \to V'$, \item natural
transformations $\overline{\phi}_{X, Y}: \phi(X) \otimes' \phi(Y)
\to \phi(X \otimes Y)$, and \item a morphism $\phi_0 : I' \to
\phi(I)$
\end{itemize}
satisfying the usual coherence conditions (see, for example,
\cite{M}). A \emph{strong} (resp. \emph{strict}) monoidal functor is
a monoidal functor in which $\overline{\phi}$ and $\phi_0$ are both
isomorphisms (resp. identities).

Recall finally that a monoidal natural transformation $$\tau:
\Phi=(\phi, \overline{\phi}, \phi_0) \to \Phi'=(\phi',
\overline{\phi'}, \phi'_0): \V \to \V'$$ between monoidal functors
is an ordinary natural transformation $\tau:\phi \to \phi'$ such
that the diagrams
$$
\xymatrix{ \phi(X) \otimes' \phi(Y) \ar[d]_{\tau_X \otimes \tau_Y}
\ar[r]^-{\overline{\phi}_{X,Y}}& \phi (X \otimes Y) \ar[d]^{\tau_{X \otimes Y}}\\
\phi'(X) \otimes' \phi'(Y) \ar[r]_-{\overline{\phi'}_{X,Y}}& \phi'
(X \otimes Y) }
$$ and $$
\xymatrix{ & \phi(I) \ar[dd]^{\tau_I}\\
I' \ar[ru]^{\phi_0} \ar[rd]_{\phi'_0}&\\
& \phi'(I)}$$ commute.

Note that the two conditions that an icon $\tau :\Phi, \Phi' :
\mathbb{B} \to \mathbb{B'}$ satisfies guarantee that
$$\tau_{\A, \A}: \phi_{\A, \A} \to \phi'_{\A, \A}: \mathbb{B}(\A,
\A)\to \mathbb{B'}(\phi(\A), \phi(\A))$$ is a monoidal natural
transformation.

Given a monoidal category $\V=(V, \otimes, I)$, we write $\V^t$ for
the monoidal category $V, \otimes^t, I)$, in which $X \otimes^t Y =Y
\otimes X.$ It is clear that if $\mathbb{B}$ is an arbitrary
bicategory, then for each $\A \in \text{Ob}(\mathbb{B})$,
$\mathbb{B}^t(\A, \A)=(\mathbb{B}(\A, \A))^t$.

\medskip

Fix a monoidal category $\V =(V, \otimes, I, \alpha, l, r)$. Recall
that a monoid in $\V$ (or $\V$-monoid) is an object $S$ of $V$
equipped with a multiplication $m_{S}: S \otimes S \to S$ and a unit
$e_{S} : I \to S$ subject to the condition that the following
diagrams commute:
$$
\xymatrix{ S\otimes (S \otimes S)\ar[d]_{S \otimes m_S}
\ar[r]^{\alpha^{-1}_{S,S,S}}& (S\otimes S) \otimes S  \ar[rr]^{m_S
\otimes S}&& S \otimes S
\ar[d]^{m_S} \\
S \otimes S \ar[rrr]_{m_S} &&& S,}
$$

$$
\xymatrix{   I \otimes S \ar[rrd]_{l_S}\ar[rr]^{e_S \otimes S} && S
\otimes S \ar[d]_{m_S}&& S \otimes I \ar[ll]_{S \otimes e_S}
\ar[lld]^{r_S} \\
 &&S && }
$$

Dually, a comonoid in $\V$ (or $\V$-comonoid) is an object
$\mathfrak{C}$ of $V$ equipped with a comultiplication
$\delta_{\mathfrak{C}}:
 \mathfrak{C} \to \mathfrak{C} \otimes \mathfrak{C}$ and a counit
$\varepsilon_{\mathfrak{C}} : \mathfrak{C} \to I$ subject to the
condition that the following diagrams commute:
$$
\xymatrix{ \mathfrak{C}\otimes (\mathfrak{C} \otimes \mathfrak{C}) &
(\mathfrak{C}\otimes
\mathfrak{C})\ar[l]_{\alpha_{\mathfrak{C},\mathfrak{C},\mathfrak{C}}}
\otimes \mathfrak{C} && \mathfrak{C} \otimes \mathfrak{C}
\ar[ll]_-{\delta_\mathfrak{C} \otimes \mathfrak{C}}
 \\
\mathfrak{C} \otimes \mathfrak{C} \ar[u]^{\mathfrak{C} \otimes
\delta_\mathfrak{C}} &&& \mathfrak{C} \ar[lll]^{\delta_\mathfrak{C}}
\ar[u]_{\delta_\mathfrak{C}},}
$$

$$
\xymatrix{  I \otimes \mathfrak{C} \ar[rrd]_{l_\mathfrak{C}} &&
\mathfrak{C} \otimes \mathfrak{C} \ar[ll]_{\varepsilon_\mathfrak{C}
\otimes \mathfrak{C}}\ar[rr]^-{\mathfrak{C} \otimes
\varepsilon_\mathfrak{C}} && \mathfrak{C} \otimes I
\ar[lld]^{r_\mathfrak{C}} \\
 &&\mathfrak{C} \ar[u]_{\delta_\mathfrak{C}} && }
$$

A morphism from a $\V$-monoid $\mathbb{S}=(S, e_S, m_S)$ to a
$\V$-monoid $\mathbb{S}'=(S', e_{S'}, m_{S'})$ is a morphism $f: S
\to S'$ in $V$ such that the diagrams
$$
\xymatrix{ I \ar[rd]_{e_{S'}}\ar[r]^{e_S} &S \ar[d]^{f}\\
& S'}$$ and $$ \xymatrix{ S \otimes S \ar[r]^{f \otimes f}
\ar[d]_{m_S}& S' \otimes S' \ar[d]^{m_{S'}}\\
S \ar[r]_{f} & S'}$$ commute.

Dually, a morphism from a $\V$-comonoid $\mathfrak{C}=(\mathfrak{C},
\varepsilon_{\mathfrak{C}}, \delta_{\mathfrak{C}})$ to a
$\V$-comonoid $\mathfrak{C}'=(\mathfrak{C}',
\varepsilon_{\mathfrak{C}'}, \delta_{\mathfrak{C}'})$ is a morphism
$\tau: \mathfrak{C} \to \mathfrak{C}'$ in $V$ such that the diagrams
$$
\xymatrix{ \mathfrak{C} \ar[r]^{\varepsilon_{\mathfrak{C}}}
\ar[rd]_{\varepsilon_{\mathfrak{C}'}}&\mathfrak{C} \ar[d]^{\tau}\\&
\mathfrak{C}}$$ and
$$ \xymatrix{ \mathfrak{C} \otimes \mathfrak{C} \ar[r]^{\tau \otimes \tau}
& \mathfrak{C}' \otimes \mathfrak{C}' \\
\mathfrak{C} \ar[r]_{\tau} \ar[u]^{\delta_{\mathfrak{C}}}&
\mathfrak{C}' \ar[u]_{\delta_{\mathfrak{C}'}}}$$ commute.

We write $\text{\textbf{Mon}}(\V)$ (resp.
$\text{\textbf{Comon}}(\V)$) for the category of $\V$-(co)monoids
and their morphisms.

If $\mathbb{S}=(S, e_S, m_S)$ is a monoid  in a monoidal category
$\V$, we write $\textbf{T}^l_{\mathbb{S}}$ (resp.
$\textbf{T}^r_{\mathbb{S}}$) for the monad on $V$ defined by
\begin{itemize}
\item $\textbf{T}^l_{\mathbb{S}}(X)=S \otimes X$ for all $X \in
\V$, \item $\eta^l_X=(e_S \otimes X) \cdot (l_X)^{-1}: X \to S
\otimes X =\textbf{T}^l_{\mathbb{S}}(X)$, and \item $\mu^l_X =(m_S
\otimes X) \cdot
(\alpha_{S,S,X})^{-1}:\textbf{T}^l_{\mathbb{S}}(\textbf{T}^l_{\mathbb{S}}(X))=S\otimes
(S\otimes X) \to S \otimes X=\textbf{T}^l_{\mathbb{S}}(X)$
\end{itemize} (resp. \begin{itemize}
\item $\textbf{T}^r_{\mathbb{S}}(X)=X \otimes S$ , \item
$\eta^r_X=(X \otimes e_S) \cdot (r_X)^{-1}: X \to X \otimes S
=\textbf{T}^r_{\mathbb{S}}(X)$, and \item $\mu_X =(X \otimes
m_S)\cdot \alpha_{X, S,
S}:\textbf{T}^r_{\mathbb{S}}(\textbf{T}^r_{\mathbb{S}}(X))=X\otimes
S\otimes S \to X \otimes S =\textbf{T}^l_{\mathbb{S}}(X))$
\end{itemize} for all
$X \in \V$.

Dually, for a $\V$-comonoid $\mathfrak{C}=(\mathfrak{C},
\varepsilon_{\mathfrak{C}}, \delta_{\mathfrak{C}})$, one defines the
comonads $\textbf{G}^l_{\mathfrak{C}}$ and
$\textbf{G}^r_{\mathfrak{C}}$ on $\V$.

\medskip

We shall need the following propositions whose proofs are
straightforward.

\begin{proposition} Any monoidal functor $$\Phi=(\phi, \overline{\phi}, \phi_0):\V \to
\V'$$ preserves monoids in the sense that, if  $\mathbb{S} = (S,
e_S, m_S)$ is a $\V$-monoid, then the triple
$\phi(\mathbb{S})=(\phi(S), \phi(e_S)\cdot \phi_0, \phi(m_S)\cdot
\overline{\phi}_{S,S})$ is a $\V'$-monoid. Moreover, if $\Phi$ is a
strong monoidal functor, then it also preserves comonoids in the
sense that, if $\mathfrak{C}=(\mathfrak{C},
\varepsilon_\mathfrak{C}, \delta_\mathfrak{C})$ is a $\V$-comonoid,
then the triple $\phi(\mathfrak{C})=(\phi(\mathfrak{C}),\phi_0^{-1}
\cdot \phi(\varepsilon_\mathfrak{C}),
(\overline{\phi}_{\mathfrak{C}, \mathfrak{C}})^{-1} \cdot
\phi(\delta_\mathfrak{C}))$ is a $\V'$-comonoid.
\end{proposition}

\begin{proposition}\label{isoiso}
1) Let $\tau: \mathfrak{C}=(\mathfrak{C},
\varepsilon_{\mathfrak{C}},
\delta_{\mathfrak{C}})\to\mathfrak{C}'=(\mathfrak{C}',
\varepsilon_{\mathfrak{C}'}, \delta_{\mathfrak{C}'})$ be an
isomorphism of $\V$-comonoids. Then the rule $$(\sigma: \mathfrak{C}
\to \mathfrak{C})\longrightarrow (\tau \sigma \tau^{-1}:
\mathfrak{C}' \to \mathfrak{C}')$$ defines an isomorphism of monoids
$$t^\tau_\V:\text{\textbf{Comon}}(\V)(\mathfrak{C},\mathfrak{C})\to \text{\textbf{Comon}}
(\V)(\mathfrak{C}',\mathfrak{C}').$$

2) Let $$\Phi=(\phi, \overline{\phi}, \phi_0):\V \to \V'$$ be a
strong monoidal equivalence. Then for any $\V$-comonoid
$\mathfrak{C}=(\mathfrak{C}, \varepsilon_{\mathfrak{C}},
\delta_{\mathfrak{C}})$, the rule
$$(\sigma: \mathfrak{C} \to
\mathfrak{C})\longrightarrow
(\phi_{\mathfrak{C},\mathfrak{C}}(\sigma): \phi(\mathfrak{C}) \to
\phi(\mathfrak{C}))$$ defines an isomorphism of monoids
$$s^\Phi_\mathfrak{C}:\text{\textbf{Comon}}(\V)(\mathfrak{C},\mathfrak{C})\to \text{\textbf{Comon}}
(\V')(\phi(\mathfrak{C}),\phi(\mathfrak{C})).$$
\end{proposition}

Given a monoid $\mathbb{S}=(S, e_S, m_S)$ in a monoidal category
$\V=(V, \otimes, I, \alpha, l, r)$, we write $\textbf{I}^l_\V
(\mathbb{S})$ for the subclass of $\text{Sub}_{\V}(S)$ consisting of
those elements $[(J, i_J : J \to S)]\in \text{Sub}_{\V}(S)$ for
which the composite
$$
\xymatrix{ \xi^l_{i_J}: S \otimes J \ar[r]^-{S \otimes i_J}& S
\otimes S \ar[r]^-{m_S}& S }$$ is an isomorphism.

Dually, we let $\textbf{I}^r_\V (\mathbb{S})$ denote the subclass of
$\text{Sub}_{\V}(S)$ consisting of those subobjects $[(J, i_J : J
\to S)]$ of $S$ for which the composite
$$
\xymatrix{ \xi^r_{i_J}: J \otimes S \ar[r]^-{i_J \otimes S}& S
\otimes S \ar[r]^-{m_S}& S }$$ is an isomorphism.

\begin{proposition}Let $f:\mathbb{S}=(S, e_S, m_S) \to \mathbb{S}'=(S', e_{S'},
m_{S'})$ be an isomorphism of $\V$-monoids. Then the assignment
$$[{J, i_J : J \to S}] \longrightarrow [(J, fi_J : J \to S')]$$
yields  bijections $$\mathbf{I}^l_\V (f): \mathbf{I}^l_\V
(\mathbb{S}) \to \mathbf{I}^l_\V (\mathbb{S}')$$ and
$$\mathbf{I}^r_\V (f): \mathbf{I}^r_\V (\mathbb{S}) \to \mathbf{I}^r_\V (\mathbb{S}').$$
\end{proposition}
\begin{proof} By symmetry, it is enough to prove the first statement.

It is clear that $\mathbf{I}^r_\V (f)$ is injective. To show that it
is surjective, consider any $[(J, i_J: J \to S')]\in \mathbf{I}^r_\V
(\mathbb{S}')$. The following diagram
$$
\xymatrix{ S \w J \ar[ddr]_{\xi^l_{i_J}}\ar[r]^{S \w i_J}& S \w S'
\ar[d]^{f \w S'}
\ar[r]^{S \w f^{-1}}& S \w S \ar[dd]^{m_S}\\
& S' \w S' \ar[d]^{m_{S'}}&&\\
&S' \ar[r]_{f^{-1}}& S}$$ commutes since $f$ (and hence also
$f^{-1}$) is a morphism of $\V$-monoids. It follows that $$m_S \cdot
(S \w f^{-1})\cdot (S \w i_J)=m_S \cdot (S \w (f^{-1}i_J))$$ is an
isomorphism. Since $f^{-1}i_J$ is clearly injective, it follows that
$[(J, f^{-1}i_J: J \to S)]\in \mathbf{I}^r_\V (\mathbb{S})$. Quite
obviously $\mathbf{I}^r_\V (f)([(J, f^{-1}i_J)])=[(J, i_J: J \to
S')]$. This completes the proof.
\end{proof}

\bigskip

The following Proposition, whose proof is left to the reader,
connects the definitions of $\mathbf{I}^l(\mathbb{S})$ to the sets
of subobjects studied in Section \ref{subcoalgebra}.

\begin{proposition}For an arbitrary monoid $\mathbb{S}=(S, e_S,
m_S)$ in a monoidal category $\V$, the assignments $$[(J, i_J : J
\to S \w I)]\longrightarrow [(J, r_S i_J : J \to S)]$$ and $$[(J,
i_J : J \to S \w I)]\longrightarrow [(J, l_S i_J : J \to S)]$$ yield
bijections
$$q^l_{\mathbb{S}}:\mathbf{Sub}_{\, \V, \,
\mathbf{T}^{\,\,l}_{\mathbb{S}}} (\mathbf{T}^{\,l}_{\mathbb{S}}(I))
\to \textbf{\emph{I}}^l_{\V}(\mathbb{S})$$ and
$$q^r_{\mathbb{S}}:\mathbf{Sub}_{\, \V, \, \textbf{T}^{\,\,r}_{\mathbb{S}}}
(\mathbf{T}^{\,r}_{\mathbb{S}}(I)) \to
\textbf{\emph{I}}^r_{\V}(\mathbb{S}),$$ respectively. When $\V$ is
strict, then $q^l_{\mathbb{S}}$ and $q^r_{\mathbb{S}}$ are
identities.
\end{proposition}

Recall that a morphism $i: X \to Y $ in a monoidal category $\V=(V,
\otimes, I)$ is said to be \emph{left} (resp. \emph{right})
\emph{pure} if for any $Z \in \V$, the morphism $Z \otimes i :Z
\otimes X \to Z \otimes Y $ (resp. $i \otimes Z : X \otimes Z \to Y
\otimes Z)$ is a monomorphism. When the unit of a monoid
$\mathbb{S}$ in $\V$ is left pure, $\mathbf{I}^l_{\V}(\mathbb{S})$
is a monoid in the sense that it has an associative product and a
unit, as the following fundamental Proposition shows.

\begin{proposition}\label{monoide} Let $\mathbb{S}=(S, e_S,
m_S)$ be a monoid in a monoidal category $\V$ such that the morphism
$e_S : I \to S$ is left pure. Then
$\textbf{\emph{I}}^l_{\V}(\mathbb{S})$ has the structure of a monoid
where the product $[(J_1, i_{J_1})]\cdot [(J_2, i_{J_12})]$ of two
elements $[(J_1, i_{J_1})], [(J_2, i_{J_12})] \in
\textbf{\emph{I}}^l_{\V}(\mathbb{S})$ is defined to be the
equivalence class of the pair $(J_1 \otimes J_2, i_{J_1 \otimes
J_2})$, where $i_{J_1 \otimes J_2}$ is the composite
$$\xymatrix{J_1 \otimes J_2 \ar[rr]^-{i_{J_1}\otimes i_{J_2}}&& S \otimes S \ar[r]^-{m_S}& S}.
$$Moreover, the element $[(I, e_S : I \to S)]$ is a two sided unit
for this multiplication.
\end{proposition}

\begin{proof} Let $[(J_1, i_{J_1})], [(J_2, i_{J_12})] \in
\textbf{\emph{I}}^l_{\V}(\mathbb{S})$. Consider the following
diagram
$$
\xymatrix{ S \otimes( J_1 \otimes J_2 )\ar[rr]^-{S \otimes (i_{J_1}
\otimes J_2)} \ar[dd]_{(\alpha_{S, J_1, J_2})^{-1}}&& S \otimes (S
\otimes J_2) \ar[rr]^-{S\otimes (S \otimes i_{J_2})}
\ar[dd]_{(\alpha_{S, S, J_2})^{-1}}&& S\otimes (S \otimes S)
\ar[rr]^-{S \otimes m_S} \ar[dd]_{(\alpha_{S, S, S})^{-1}} && S
\otimes S
\ar[dddd]^{m_S}\\\\
(S \w J_1) \w J_2 \ar@{}[rruu]^{(1)}\ar[rr]_{(S \w i_{J_1}) \w J_2}
&& (S\w S)\w J_2 \ar@{}[rruu]^{(2)} \ar[rr]_{(S\w S)\w i_{J_2}}
\ar[dd]_{m_S \otimes J_2}&& (S \w S)
\w S  \ar[dd]_{m_S \otimes S} \ar@{}[rr]^{(4)} &&&\\\\
&& S \otimes J_2 \ar@{}[rruu]^{(3)}\ar[rr]_-{S \otimes i_{J_2}} && S
\otimes S \ar[rr]_-{m_S} && S}
$$ in which
\begin{itemize}
\item (1) and (2) are commutative by naturality of $\alpha$;

\item (3) commutes by naturality of $m_S \otimes -$, and

\item (4) commutes by associativity of $m_S$.
\end{itemize}

It follows that
$$\xi_{i_{J_1 \otimes J_2}}=m_S \cdot (S \w i_{J_1 \w J_2})=
m_S \cdot (S \w m_S)\cdot (S \w (i_{J_1}\w i_{J_2}))=$$
$$=m_S \cdot (S \otimes m_S)\cdot (S
\otimes (S\otimes i_{J_2}))\cdot (S \otimes (i_{J_1} \otimes
J_2))=$$$$=m_S \cdot (S \otimes i_{J_2}) \cdot (m_S \otimes
J_2)\cdot ((S \otimes i_{J_1})\otimes J_2)\cdot (\alpha_{S, J_1,
J_2})^{-1} =$$$$=\xi_{i_{J_2}}\cdot ((m_S \cdot (S \otimes
i_{J_1}))\otimes J_2)\cdot (\alpha_{S, J_1,
J_2})^{-1}=\xi_{i_{J_2}}\cdot (\xi_{i_{J_1}}\otimes
J_2)\cdot(\alpha_{S, J_1, J_2})^{-1}.$$ But the morphisms
$\xi_{i_{J_1}}$ and $\xi_{i_{J_2}}$ are both isomorphisms, since
$[(J_1, i_{J_1})], [(J_2, i_{J_12})] \in
\mathbf{I}^l_{\V}(\mathbb{S})$, implying that $\xi_{i_{J_1 \otimes
J_2}}$ is also an isomorphism.

Next, since $\xi_{i_{J_1 \otimes J_2}}=m_S \cdot (S \otimes i_{J_1
\otimes J_2})$ is an isomorphism, the morphism $S \otimes i_{J_1
\otimes J_2}$ is a monomorphism, and since $e_S : I \to S$ is left
pure by our assumption, the functor $S \otimes - : V \to V$ reflects
monomorphisms, implying  that  $i_{J_1 \otimes J_2}$ is a
monomorphism. Thus $[(J_1 \otimes J_2, i_{J_1 \otimes J_2})] \in
\mathbf{I}^l_{\V}(\mathbb{S}).$

Suppose now that $[(J_1, i_{J_1})]=[(J'_1, i_{J'_1})]$ and $[(J_2,
i_{J_2})]=[(J'_2, i_{J'_2})]$. Then there exist isomorphisms $f_1:
J_1\to J'_1$ and $f_2: J_2\to J'_2$ such that the diagrams

\begin{equation}
\xymatrix{ J_1 \ar[dd]_{f_1} \ar[rrd]^{i_{J_1}} && \\
&&  S\\
J'_1 \ar[rru]_{i_{J'_1}}&&}
\end{equation} and

\begin{equation}
\xymatrix{ J_2 \ar[dd]_{f_2} \ar[rrd]^{i_{J_2}} && \\
&&  S\\
J'_2 \ar[rru]_{i_{J'_2}}&&} \end{equation} commute. Then $$ i_{J'_1
\otimes J'_2}
 \cdot (f_1 \otimes f_2)=m_S \cdot (S \otimes i_{J'_2} )\cdot
 (i_{J'_1}\otimes J'_2)\cdot (f_1 \otimes J'_2)\cdot (J_1 \otimes f_2)=\,\,{\text{by} \,\,\, (6)}$$
 $$= m_S \cdot (S \otimes i_{J'_2} )\cdot
 (i_{J_1}\otimes J'_2)\cdot (J_1 \otimes f_2)=\,\,{\text{by \,\,
 naturality}}$$
 $$=m_S \cdot (S \otimes i_{J'_2} )\cdot
 (S \otimes f_2)\cdot (i_{J_1} \otimes J_2)= \,\,{\text{by} \,\,\,
(7)}$$$$= m_S \cdot (S \otimes i_{J_2} )\cdot
  (i_{J_1} \otimes J_2)=i_{J_1
\otimes J_2}$$ and hence the following diagram

$$
\xymatrix{ J_1 \otimes J_2 \ar[dd]_{f_1 \otimes f_2}\ar[drr]^{i_{J_1
\otimes J_2}} && \\
&& S \\
 J'_1 \otimes J'_2 \ar[rru]_{i_{ J'_1 \otimes J'_2}}, && }
$$ commutes, whose commutativity just means -- since $f_1 \otimes
f_2$ is an isomorphism -- that $[(J_1 \otimes J_2, i_{J_1 \otimes
J_2})]=[(J'_1 \otimes J'_2, i_{J'_1 \otimes J'_2})]$. Thus, the map
$([(J_1, i_{J_1})], [(J_2, i_{J_2})]) \to [(J_1 \otimes J_2, i_{J_1
\otimes J_2})]$ is well-defined.

Since $e_S : I \to S$ is left pure, $e_S$ is a monomorphism.
Moreover, since $e_S$ is the unit for the multiplication $m_S : S
\otimes S \to S$, the following diagram
$$
\xymatrix{ S \otimes I \ar[rd]_{r_S}\ar[r]^{S \otimes e_S} & S
\otimes S \ar[d]^{m_S}\\
&S} $$ commutes. And since $r_S$ is an isomorphism, $[(I, e_S)] \in
\mathbf{I}^l_{\V}(\mathbb{S})$.

Now, for any $[{J, i_J}] \in \mathbf{I}^l_{\V}(\mathbb{S})$,
considering the diagram
$$\xymatrix{
I \otimes J  \ar[rr]^{I \otimes i_J} \ar[dd]_{l_J} && I \otimes S
\ar[rr]^{e_S \otimes S} \ar[dd]_{l_S}&& S \otimes S
\ar[ddll]^{m_S}\\\\
J \ar[rr]_{i_J}&&S && }$$ (resp.
$$\xymatrix{
J \otimes I  \ar[rr]^{i_J \otimes I} \ar[dd]_{r_J} && S \otimes I
\ar[rr]^{S\otimes e_S } \ar[dd]_{r_S}&& S \otimes S
\ar[ddll]^{m_S}\\\\
J \ar[rr]_{i_J}&&S && })$$ which is  commutative, since

\begin{itemize}
\item the square commutes by naturality of $l$ (resp. $r$), and
\item the triangle commutes because $e_S$ is the unit for $m_S : S
\otimes S \to S$,
\end{itemize} one sees that $[{I, e_S}]$ is a two sided unit for
the multiplication
$$([(J_1, i_{J_1})], [(J_2, i_{J_2})]) \to [(J_1 \otimes J_2,
i_{J_1 \otimes J_2})].$$

\end{proof}

\begin{remark} Let $\mathbb{S}=(S, e_S,
m_S)$ be an arbitrary $\V$-monoid. If $e_S : I \to S$ is left (resp.
right) pure, then transporting the structure of a monoid on
$\textbf{{I}}^l_{\V}(\mathbb{S})$ (resp.
$\textbf{{I}}^r_{\V}(\mathbb{S})$) to $\text{\textbf{Sub}}_{\, \V,
\, \textbf{T}^{\,\,l}_{\mathbb{S}}}
(\textbf{T}^{\,l}_{\mathbb{S}}(I))$ (resp. $\text{\textbf{Sub}}_{\,
\V, \, \textbf{T}^{\,\,l}_{\mathbb{S}}}
(\textbf{T}^{\,r}_{\mathbb{S}}(I)$), one sees that
$q^l_{\mathbb{S}}$ and $q^r_{\mathbb{S}}$ become isomorphisms of
monoids.
\end{remark}

\begin{proposition}\label{funtormorfismo}
 Let $$\Phi=(\phi, \overline{\phi}, \phi_0): \V \to \V'$$
be a strong monoidal functor such that the functor $\phi : V \to V'$
preserves monomorphisms. Then for any monoid $\mathbb{S}=(S, e_S,
m_S)$ in $\V$, the assignment $$[(J, i_J)] \longrightarrow
[(\phi(J), \phi(i_J))]$$ yields a map
$$\phi_{\,\mathbb{S}}:\textbf{\emph{I}}^l_{\V}({\mathbb{S}}) \to
\textbf{\emph{I}}^l_{\V'}({\phi(\mathbb{S}))}.$$
\end{proposition}
\begin{proof} Given an arbitrary $[(J, i_J)]\in
\textbf{{I}}^l_{\V}(\mathbb{S})$, consider the pair $(\phi(J),
\phi(i_J))$. Since the functor $\phi$ preserves monomorphisms,
$\phi(i_J): \phi(J) \to \phi(S)$ is a monomorphism and thus
$[(\phi(J), \phi(i_J))] \in \text{Sub}_{\V'}(\phi(S))$. To see that
$[(\phi(J), \phi(i_J))] \in \textbf{{I}}^l_{\V'}(\mathbb{\phi(S)})$,
consider the composite
$$
\xymatrix{ \xi_{\phi(i_J)}: \phi(S)\otimes \phi(J)
\ar[rr]^-{\phi(S)\otimes \phi(i_J)}&& \phi(S)\otimes \phi(S)
\ar[rr]^-{m_{\phi(S)}}&&\phi(S).}$$ Since $m_{\phi(S)}=\phi(m_S)
\cdot \overline{\phi}_{S,S}$, we have $$\xi_{\phi(i_J)}=\phi(m_S)
\cdot \overline{\phi}_{S,S} \cdot (\phi(S)\otimes
\phi(i_J))=\,\,\text{by \,\, naturality \,\,of \,\,}
\overline{\phi}$$$$=\phi(m_S) \cdot \phi(S \otimes i_J)\cdot
\overline{\phi}_{S, J}=$$$$=\phi(m_S \cdot (S \otimes i_J))\cdot
\overline{\phi}_{S, J}=$$$$=\phi(\xi_{i_J})\cdot \overline{\phi}_{S,
J}.$$ Thus $\xi_{\phi(i_J)}$ is an isomorphism, implying that
$[(\phi(J), \phi(i_J))] \in
\textbf{{I}}^l_{\V'}({\phi(\mathbb{S})})$.
\end{proof}

\begin{proposition}\label{reduccionaestricto}
In the situation of Proposition \ref{funtormorfismo}, if $\Phi$ is a
strong monoidal equivalence, then the map $\phi_{\,\mathbb{S}}$ is
bijective. Moreover, when $e_S : I \to S$ is left pure, this
bijection is an isomorphism of monoids.
\end{proposition}
\begin{proof} Quite obviously, the map $\phi_{\,\mathbb{S}}$ is
injective. So it suffices to show that $\phi_{\,\mathbb{S}}$ is
surjective. So suppose that $[(X,i_X)]\in
\textbf{{I}}^l_{\V'}({\phi(\mathbb{S})})$, and let $\Phi'=(\phi',
\overline{\phi'}, \phi'_0): \V' \to \V$ be a strong monoidal functor
with monoidal natural isomorphisms $\sigma: 1 \to \phi\phi'$ and
$\tau: \phi'\phi \to 1$. Since any equivalence preserves
monomorphisms, it follows from Proposition \ref{funtormorfismo} that
$[(\phi'(i_X): \phi'(X) \to \phi'\phi(S))]\in
\textbf{{I}}^l_{\V}(\phi'\phi(\mathbb{S})).$ It is then clear that
$(\tau_S \cdot \phi'(i_X): \phi'(X)\to S) \in \text{Sub}_\V (S).$
Next, looking at the following diagram
$$
\xymatrix{ S \otimes \phi'(X) \ar[rr]^-{S \otimes \phi'(i_X)}
\ar[dd]_{(\tau_{S})^{-1}\otimes \phi'(X)} && S \otimes \phi'\phi(X)
\ar[rr]^-{S \otimes \tau_S} \ar[dd]^{(\tau_{S})^{-1}\otimes
\phi'\phi(S)}&& S \otimes S
\ar[dd]^{m_S}\\\\
\phi'\phi(S) \ar@{}[rruu]_{(1)}\ar[dd]_{\overline{\phi'}_{\phi(S),
X}}\ar[rr]_-{\phi'\phi(S) \otimes \phi'\phi(i_X)}&& \phi'\phi(S)
\otimes \phi'\phi(S) \ar[dd]^{\overline{\phi'}_{\phi(S),
\phi(S)}}\ar@{}[rr]_{(2)}&& S \\\\
\phi'(\phi(S))\otimes X
\ar@{}[rruu]_{(3)}\ar[rr]_{\phi'(\phi(S))\otimes i_X}&&
\phi'(\phi(S) \otimes \phi(S)) \ar[rr]_-{\phi'(m_S)}&& \phi'\phi(S)
\ar[uu]_{\tau_S}}
$$ which is commutative since
\begin{itemize}
\item (1) commutes by naturality, \item (2) commutes, since
$\tau_S : \phi'\phi(S) \to S$ is a morphism of monoids (which
follows from the fact that $\tau$ is a monoidal natural
transformation), and \item (3) commutes since $\tau$ is a monoidal
natural transformation,
\end{itemize}one sees that $\xi_{i_{\tau_S \cdot \phi'(i_X)}}=m_S \cdot
(S \otimes \tau_S)\cdot (S \otimes \phi'(i_X))$ is an isomorphism,
proving that $[(\tau_S \cdot \phi'(i_X): \phi'(X)\to S)] \in
\textbf{{I}}^l_{\V}({\mathbb{S}})$.

We now have:
$$\phi(\tau_S \cdot \phi'(i_X))=\phi(\tau_S) \cdot
\phi(\phi'(i_X))=\,\,\,\,\text{since}\,\, \phi\tau \cdot\sigma
\phi=1$$$$=(\sigma)^{-1}_{\phi(S)}\cdot
\phi(\phi'(i_X))=\,\,\text{by \,naturality\, of}\,\, \sigma$$$$=i_X
\cdot (\sigma)^{-1}_X . $$ Thus the diagram
$$
\xymatrix{ \phi\phi'(X)\ar[dd]_{(\sigma)^{-1}_X }
\ar[rrd]^-{\phi(\tau_S \cdot
\phi'(i_X))}&&\\
&& \phi(S)\\
X \ar[rru]_{i_X}&&}$$ commutes, proving that
$\phi_{\,\mathbb{S}}([(\phi'(X),\tau_S \cdot \phi'(i_X) )])=[(X,
i_X)]$. Thus the map
$$\phi_{\,\mathbb{S}}:\textbf{{I}}^l_{\V}({\mathbb{S}}) \to
\textbf{{I}}^l_{\V'}({\phi(\mathbb{S}))}$$ is surjective, and hence
bijective.

Suppose now that the morphism $e_S :I \to S$ is left pure in $\V$.
It is easy to see- using the fact that $\Phi$ is a strong monoidal
equivalence- that the morphism $\phi(e_S): \phi(I) \to \phi(S)$ is
left pure in $\V'$. Since $\phi_0 : I' \to \phi(I)$ is an
isomorphism, this implies that $\phi(e_S)\cdot \phi_0 : I' \to
\phi(S)$ is also left pure in $\V'$. Consequently, when the morphism
$e_S :I \to S$ is left pure in $\V$, then
$\textbf{{I}}^l_{\V}({\mathbb{S}})$ and
$\textbf{{I}}^l_{\V'}({\phi(\mathbb{S}))}$ both have the structure
of a monoid. We want to show that $\phi_{\,\mathbb{S}}$ is a
homomorphism (and hence an isomorphism) of monoids.

Quite obviously, the diagram
$$
\xymatrix{ I' \ar[dd]_{\phi_0} \ar[rrd]^-{\phi(e_S)\cdot \phi_0}&&\\
&& \phi(S)\\
\phi(I) \ar[rru]_{\phi(e_S)}&&}$$ commutes, proving that
$\phi_{\,\mathbb{S}}([(I, e_S)])=[(I', e_{\phi(S)})]$. Thus the map
$\phi_{\,\mathbb{S}}$ preserves the identity.

Now let $[(J_1, i_{J_1})], [(J_2, i_{J_2})] \in
\textbf{{I}}^l_{\V}({\mathbb{S}})$ be arbitrary elements.  Since in
$\textbf{{I}}^l_{\V}({\mathbb{S}})$,
$$[(J_1, i_{J_1})]\cdot[(J_2, i_{J_2})]=[(J_1 \otimes J_2, i_{J_1
\otimes J_2})],$$ where
$$i_{J_1 \otimes J_2}=m_S \cdot (i_{J_1}\otimes i_{J_2}),$$ we
have $$\phi_{\mathbb{S}}([(J_1, i_{J_1})]\cdot[(J_2,
i_{J_2})])=[(\phi(J_1 \otimes J_2), \phi(m_S \cdot (i_{J_1}\otimes
i_{J_2}))]=$$$$=[(\phi(J_1 \otimes J_2), \phi(m_S) \cdot
\phi(i_{J_1}\otimes i_{J_2}))].$$ Considering the diagram
$$
\xymatrix{ \phi(J_1)\otimes \phi(J_2) \ar[rr]^{\phi(i_{J_1})\otimes
\phi(i_{J_2})} \ar[dd]_{\overline{\phi}_{J_1,\, J_2}}&& \phi(S)
\otimes \phi(S) \ar[dd]^{\overline{\phi}_{S, S}} \\\\
\phi(J_1 \otimes J_2) \ar[rr]_{\phi(i_{J_1} \otimes i_{J_2})}&&
\phi(S \otimes S)}$$ which is commutative by naturality of
$\overline{\phi}$, one sees that $$ \phi(i_{J_1 \otimes J_2})\cdot
\overline{\phi}_{J_1,\, J_2} =\phi(m_S) \cdot \phi(i_{J_1}\otimes
i_{J_2})\cdot \overline{\phi}_{J_1,\, J_2}=$$$$=\phi(m_S) \cdot
\overline{\phi}_{S,S}\cdot (\phi(i_{J_1})\otimes \phi(i_{J_2}))=
$$$$=m_{\phi(S)}\cdot (\phi(i_{J_1})\otimes \phi(i_{J_2}))=i_{\phi(J_1)\otimes
\phi(J_2)},$$which means that the diagram
$$
\xymatrix{\phi(J_1) \otimes \phi(J_2)
\ar[dd]_{\overline{\phi}_{J_1,\, J_2}}
\ar[rrd]^-{i_{\phi(J_1)\otimes
\phi(J_2)}}&&\\
&& \phi(S)\\
\phi(J_1 \otimes J_2) \ar[rru]_{\phi(i_{J_1 \otimes J_2})}&&}$$
commutes. Thus in $\mathbf{I}^l_{\V'}(\phi(\mathbb{S}))$
$$[({\phi(J_1 \otimes J_2),\phi(i_{J_1 \otimes
J_2})})]=[(\phi(J_1) \otimes \phi(J_2), i_{\phi(J_1)\otimes
\phi(J_2)})],$$ proving that $$\phi_{\mathbb{S}}([(J_1,
i_{J_1})]\cdot[(J_2, i_{J_2})])=\phi_{\mathbb{S}}([(J_1,
i_{J_1})])\cdot\phi_{\mathbb{S}}([(J_2, i_{J_2})]).$$ Thus
$\phi_{\mathbb{S}}$ is a homomorphism of monoids.
\end{proof}

\bigskip

Let us now consider a bicategory $\mathbb{B}$ and an arbitrary
0-cell $\A \in \text{Ob}(\mathbb{B})$. We call a (co)monoid in the
monoidal category $\mathbb{B}(\A, \A)$ an $\A$-(co)ring and write
$\A$-$\text{\textbf{Rings}}$ (resp. $\A$-$\text{\textbf{Corings}}$)
for the category of $\A$-rings (resp. $\A$-corings). When
$\mathbb{B} = \mathsf{Bim}$, the bicategory of bimodules over unital
rings, we recover the usual notions of an $A$--ring or an
$A$--coring for a given ring $A$.

Suppose that $\mathbb{S}=(S,e_S,m_S)$ be an $\A$-ring. For any
$0$-cell $\C$, $\mathbb{S}$ induces a monad
$\textbf{T}=\textbf{T}^{\C}_\mathbb{S}=(T,m_\textbf{T},
e_\textbf{T})$ on the category $\mathbb{B}(\C, \A)$ as
follows:$$T(f)=S \w f, \, (m_\textbf{T})_f=(m_S \w f)\cdot
(\alpha_{S,S, f})^{-1},\,\,\text{and}\,\, (e_\textbf{T})_f=(e_S \w
f)\cdot (l_f)^{-1} \,\, \text{for \,\,all} \,\, f \in \mathbb{B}(\C,
\A).$$ We write $\mathbb{B}(\C,\A)_{\textbf{T}^{\C}_\textbf{S}}$ for
the Eilenberg-Moore category of
$\textbf{T}^{\C}_\textbf{S}$-algebras.

Dually, given an $\A$-coring $\mathfrak{C}=(\mathfrak{C},
\e_\mathfrak{C}, \delta_\mathfrak{C})$, one defines a comonad
$\textbf{G}=\textbf{G}^{\C}_\mathfrak{C}=(G,
\delta_\textbf{G},\e_\textbf{G})$ on the category $\mathbb{B}(\C,
\A)$ by $$G(h)=\mathfrak{C} \w
h,\,(\delta_\textbf{G})_h=\alpha_{\mathfrak{C},\mathfrak{C},h}(\delta_\mathfrak{C}
\w h)\,\, \text{and}\,\, (\e_\textbf{G})_h=l_h \cdot
(\varepsilon_\mathfrak{C} \w h) \,\, \text{for \,\,all} \,\, h \in
\mathbb{B}(\C, \A).$$ We write $\mathbb{B}(\C,\A)^
{\textbf{G}^\C_\mathfrak{C}}$ for the corresponding Eilenberg-Moore
category of $\textbf{G}^\C_\mathfrak{C}$-coalgebras.

\medskip

 Recall that a 1-cell $f :\A \to \B$ admits as a right
adjoint a 1-cell $f^* : \B \to \A$ when there exist 2-cells
$\eta_{f}: 1_{\A} \to f^* \w f$ and $\e_{f}: f \w f^* \to 1_{\B}$
such that the following diagrams commute in $\mathbb{B}(\A, \B)$ and
$\mathbb{B}(\B, \A)$, respectively:
\begin{equation}\label{adj1} \xymatrix{ f \w 1_{\A} \ar[r]^{f \w \eta_f}\ar[d]_{r_f}
& f \w (f^* \w f)  \ar[r]^{(\alpha_{f, f^*, f})^{-1}}&
(f \w f^*)\w\phi \ar[d]^{\e_f \w f}\\
f \ar[rr]_{(l_f)^{-1}}&& 1_\B \w f}
\end{equation}
and
\begin{equation}\label{adj2}\xymatrix{ 1_\A \w f^*  \ar[r]^{\eta_f \w f^* }
\ar[d]_{l_{f^*}}& (f^* \w f )\w f^* \ar[r]^{\alpha_{f^*, f, f^*}}
& f^* \w (f \w f^*) \ar[d]^{ f^* \w \e_f } \\
f^* \ar[rr]_{(r_{f^*})^{-1}}&& f^* \w 1_\B \,.}
\end{equation}

We usually write $\eta_f, \e_f :f \dashv f^* : \B \to \A$ to denote
that the 1-cell $f^*$ is right adjoint to the 1-cell $f$ with unit
$\eta_f$ and counit $\e_f$.

When $\mathbb{B}$ is a 2-category, then \eqref{adj1} and
\eqref{adj2} can be rewritten as follows:
\begin{equation}\label{adj3}
\e_f f \cdot f \eta_f =1 \,\,\text{and} \,\, f^* \e_{f} \cdot
\eta_{f}f^*=1.
\end{equation}

Let $\eta_f, \e_f :f \dashv f^* : \B \to \A$ be adjunction in
$\mathbb{B}$ and let $\C$ be an arbitrary 0-cell of $\mathbb{B}$.
Since the representable $$\mathbb{B}(\C, -) : \mathbb{B}\to
\text{CAT}$$ is a homomorphism of bicategories and since any
homomorphism of bicategories preserves adjunctions, the functor
$$\mathbb{B}(\C,f)=f \w - :\mathbb{B}(\C, \A) \to \mathbb{B}(\C, \B)$$ admits as a right
adjoint the functor $$\mathbb{B}(\C,f^*)=f^* \w - : \mathbb{B}(\C,
\B) \to \mathbb{B}(\C, \A)\,. $$ The unit $\eta_f^\C$ and counit
$\e_f^\C$ of this adjunction are given by the formulas:
$$\xymatrix{(\eta_f^\C)_g : g \ar[r]^-{(l_g)^{-1}} & 1_\A \w g
\ar[r]^-{\eta_f \w g}& (f^* \w f) \w g \ar[r]^-{\alpha_{f^*, f, g}}&
f^* \w (f \w g),\,\,\, \text{for\,\,all}\,\, g \in \mathbb{B}(\C,
\A)}$$ and
$$\xymatrix{
(\e_f^\C)_h : f \w (f^* \w h)\ar[rr]^-{(\alpha_{f,f^*, h})^{-1}}&&(f
\w f^*) \w h \ar[r]^-{\e_f \w h} &1_\B \w h \ar[r]^-{r_h}& h, \,\,\,
\text{for\,\,all}\,\, h \in \mathbb{B}(\C, \B).}$$

We shall write $\textbf{T}^\C_f$ (resp. $\textbf{G}^\C_f$) for the
monad (resp. comonad) on the category $\mathbb{B}(\C, \A)$ (resp.
$\mathbb{B}(\C, \B)$) generated by the adjunction
$$\mathbb{B}(\C,f) \dashv \mathbb{B}(\C,f^*):\mathbb{B}(\C, \B)\to\mathbb{B}(\C, \A).$$

Dually, the functor $$\mathbb{B}(f,\C)=-\w f : \mathbb{B}(\B , \C)
\to \B(\A, \C)$$ is right adjoint to the functor
$$\mathbb{B}(f^*,\C)=-\w f^* : \mathbb{B}(\A, \C) \to
\mathbb{B}(\B, \C) $$ whose unit and counit are
 $$\xymatrix{
(\overline{\eta}_f^\C)_g : g \ar[r]^-{(r_g)^{-1}}& g \w 1_\A
\ar[r]^-{g \w \eta_f} & g \w (f^* \w f)\ar[rr]^-{(\alpha_{g,f^*,
f})^{-1}}&& (g\w f^*)\w f \,\,\text{for \,\,\,all\,\,} g:\A \to
\C,}$$ and
$$\xymatrix{ (\overline{\e}_f^\C)_h :(h \w f)\w f^*
\ar[rr]^-{\alpha_{h, f, f^*}}&& h \w ( f \w f^*) \ar[r]^-{h \w
\e_f}& h \w 1_\B \ar[r]^-{r_h}&
 h \,\,\text{for \,\,\,all\,\,} h:\B \to \C.}$$
From (7) and (8), it is easily verified that

\begin{proposition} For any adjunction $\eta_f, \e_f :
f \dashv f^* : \B \to \A$ in $\mathbb{B}$, the triple
$$S^{\mathbb{B}}_f=(f^* \w f, m_f,\eta_f ),$$
where $m_f$ is the composite
$$
(r_{f^*} \w f^*  )\cdot((f^* \w \e_f) \w f)\cdot( \alpha_{f^*, f,
f^*}\w f)\cdot(\alpha_{f^* \w f,f^*, f })^{-1}:(\phi^* \w \phi )\w
(f^* \w f)\to f^* \w f,
$$ is an $\A$-ring, while the triple
$$\mathfrak{C}^{\mathbb{B}}_f =(f \w f^*,
\delta_f , \e_f), $$ where $\delta_f $ is the composite
$$
(\alpha_{f^* \w f,f^*, f })\cdot ((\alpha_{f, f^*, f})^{-1}\w f^*)
\cdot((f \w \eta_f) \w f^* )\cdot((r_f)^{-1}\w f^*):f \w f^* \to (f
\w f^*)\w (f \w f^*),
$$is a $\B$-coring.
\end{proposition}

By abuse of notation, we will write $S_f$ and $\mathfrak{C}_f $ when
$\mathbb{B}$ clear from the context. Observe that the coring
$\mathfrak{C}^{\mathbb{B}}_f$ is a general form of the comatrix
coring introduced in \cite{EG1}.

\bigskip

\begin{remark}\label{monisomorphic}
It is not hard to check that the monads $\textbf{T}^\A_f$ and
$\textbf{T}^l_{S_f}$ are isomorphic.
\end{remark}

We now consider an adjunction $\eta_f, \e_f :f \dashv f^* : \B \to
\A$ in an arbitrary 2-category $\mathbb{K}$, the corresponding
$\B$-coring $\mathfrak{C}_f$ and write
$\text{End}_\B(\mathfrak{C}_f,\mathfrak{C}_f)$ for
$\B$-\text{\textbf{Corings}}$(\mathfrak{C}_f,\mathfrak{C}_f)$. Since
$(f, \eta_f f : f \to ff^* f) \in \mathbb{K}(\A,
\B)^{\textbf{G}^\A_{\mathfrak{C}_f}}$, the pair $(f, \alpha f \cdot
f \eta_f)$ is also an object of the category $\mathbb{B}(\A,
\B)^{\textbf{G}^\A_{\mathfrak{C}_f}} $ for any $\alpha \in
\text{End}_\B(\mathfrak{C}_f,\mathfrak{C}_f)$. Hence the assignment
$$\alpha \to (f, \alpha f \cdot f \eta_f)$$ yields a map
$$\gamma: \text{End}_\B(\mathfrak{C}_f,\mathfrak{C}_f) \to \textbf{G}^\A_{\mathfrak{C}_f}\text{-coalg}(f).$$

\begin{proposition}\label{endocoal}
The map $\gamma$ is bijective.
\end{proposition}
\begin{proof}For any $(f, \theta_f: f \to \mathfrak{C}^*_f(f)=ff^*f ) \in \mathbb{K}
(\A, \B)^{\textbf{G}^\A_{\mathfrak{C}_f}},$ write
$\gamma'(\theta_f)$ for the composite
$$\xymatrix{\mathfrak{C}_f =ff^* \ar[r]^-{\theta_f f^*}& f
f^* f f^* \ar[r]^-{ff^* \e_f}& ff^*=\mathfrak{C}_f}.$$ We claim that
$\gamma'(\theta_f)$ is a endomorphism of the $\B$--coring
$\mathfrak{C}_f$. Indeed, looking at the diagram
$$\xymatrix{
ff^* \ar[rr]^-{\theta_f f^*} \ar@{=}[drr]&& ff^*ff^*
\ar[rr]^-{ff^*\e_f } \ar[d]^-{f \e_f f^*}
&& ff^* \ar[d]^-{\e_f}&\\
&& ff^*\ar[rr]_-{\e_f}&& 1 }$$ in which the triangle commutes
because $(f, \theta_f) \in \mathbb{K}(\A,
\B)^{\textbf{G}^\A_{\mathfrak{C}_f}}$, while the square commutes by
functoriality of composition, we see that $\e_f \cdot
\gamma'(\theta_f)=\e_f.$ Moreover, since
\begin{equation}\label{adj4}
ff^*\e_f f \cdot \theta_f f^* f \cdot f \eta_f=ff^* \e_ff \cdot
ff^*f\eta_f \cdot \theta_f=\theta_f \end{equation} by functoriality
of composition and by \eqref{adj3}, we have

$$(\gamma'(\theta_f)\gamma'(\theta_f)) \cdot \delta_{\mathfrak{C}_f}=f
f^*\gamma'(\theta_f) \cdot \gamma'(\theta_f) ff^* \cdot
\delta_{\mathfrak{C}_f}=$$$$=ff^*ff^*\e_f  \cdot ff^* \theta_f f^*
\cdot ff^* \e_f ff^* \cdot \theta_f f^* ff^* \cdot f\eta_f
f^*=\,\,\text{by} \,\,\eqref{adj4} $$$$=ff^*ff^* \e_f  \cdot
ff^*\theta_f f^* \cdot \theta_f f^*= \,\,\text{since}\,\, (f,
\theta_f) \in \mathbb{K}(\A,
\B)^{\textbf{G}^\A_{\mathfrak{C}_f}}$$$$=ff^*ff^* \e_f \cdot f\eta_f
f^*ff^* \cdot \theta_f f^*= \,\,\text{by \, naturality \, of\,
composition}$$$$=f \eta_f f^* \cdot ff^* \e_f \cdot \theta_f
f^*=\delta_{\mathfrak{C}_f} \cdot \gamma'(\theta_f).$$

Thus $\gamma'(\theta_f)$ is an endomorphism of the $\B$-coring
$\mathfrak{C}_f$, and therefore the assignment $$(f, \theta_f) \to
\gamma'(\theta_f)$$ yields a map

$$\gamma': {\textbf{G}^\A_{\mathfrak{C}_f}}\text{-coalg}(f) \to \text{End}_\B(\mathfrak{C}_f,\mathfrak{C}_f).$$

We are now going to show that $\gamma'$ is the inverse to $\gamma$.
To show that $\gamma' \gamma=1$, consider an arbitrary element
$\alpha$ of
 $\text{End}_{\B}(\mathfrak{C}_f,\mathfrak{C}_f)$. We have:

$$(\gamma' \gamma) (\alpha)=\gamma'(f, \alpha f \cdot
f\eta_f)=$$$$= ff^* \e_f \cdot \alpha ff^* \cdot f \eta_f
f^*=\,\,\text{by \,naturality \,of \,composition}$$$$= \alpha \cdot
ff^* \e_f \cdot f \eta_f f^*=\,\,\text{by
\,\,}\eqref{adj3}$$$$=\alpha.$$ Thus $\gamma'\gamma=1.$

Now, if $(f, \theta_f)\in
{\textbf{G}^\A_{\mathfrak{C}_f}}{-coalg}(f)$, then we have:
$$(\gamma \gamma')(\theta_f)=\gamma(ff^* \e_f \cdot \theta_f f^*))=$$
$$=ff^* \e_f f \cdot \theta_f f^* f \cdot f \eta_f = \,\,
\text{by \,naturality \,of \,composition}$$
$$=ff^* \e_f f \cdot ff^* f \eta_f \cdot \theta_f=
\,\,\text{by \,\,}\eqref{adj3}$$$$=\theta_f.$$ Thus
$\gamma\gamma'=1$ and hence $\gamma$ is bijective whose inverse is
$\gamma'$.
\end{proof}

\bigskip
\bigskip

\begin{proposition}\label{monoide2}
Let $f \dashv f^* : \B \to \A$ be an adjunction in a bicategory
$\mathbb{B}$ such that the unit of the adjunction
$$\mathbb{B}(\A,f) \dashv \mathbb{B}(\A,f^*) : \mathbb{B}(\A,
\A)\to \mathbb{B}(\A, \B)$$ is a monomorphism. Then
$$\mathbf{Sub}_{\,\mathbb{B}(
\A,\,\A\,),\,\mathbb{B}(\A,f)}(\mathbb{B}(\A,f^*)(f))=\mathbf{Sub}_{\,\mathbb{B}(
\A,\,\A\,),\,\mathbb{B}(\A,f)}(f^*f)$$ has the structure of a monoid
where the product $$[(h', i_{h'})] \cdot [(h,i_{h})]$$ of
$$[(h,i_{h})], \, [(h',i_{h'})] \in
\mathbf{Sub}_{\,{\mathbb{B}(\A,\,\A\,),\,\mathbb{B}(\A,f)}} (f^*f)$$
is defined to be the (the equivalence class of) the pair $(h'\w h,
i_{h'\w h})$, where $i_{h'h}$ is the composite
$$\xymatrix{ h'\w h \ar[r]^-{i_{h'}\w h} & (f^* \w f)\w h
\ar[r]^-{(f^* \w f) \w i_{h}} & (f^*\w f) \w (f^*\w
f)\ar[rr]^-{(\alpha_{f^* \w f, f^*, f})^{-1}}&&  ((f^* \w f)\w
f^*)\w f \\\ar[r]^-{\alpha_{f^*, f, f^*}\w f}& (f^* \w(f \w f^*)) \w
f \ar[r]^-{(f^* \w \e_f) \w f} & (f^*\w 1_\A) \w f
\ar[r]^-{r_{f^*}\w f}& f^* \w f.}$$ Moreover, the element
$[(\eta_{f} : 1_{\A} \to f^*f)] $ is a two sided unit for this
multiplication.
\end{proposition}

\begin{proof}
According to Propositions \ref{mon1} and \ref{mon2} and Remark
\ref{monisomorphic}, we have the following chain of bijections

$$\text{\textbf{Sub}}_{\,\mathbb{B}(
\A,\,\A\,),\,\mathbb{B}(\A,f)}(\mathbb{B}(\A,f^*)(f))\simeq
\text{\textbf{Sub}}_{\,\mathbb{B}(
\A,\,\A\,),\,\mathbb{B}(\A,f)}(\mathbb{B}(\A,f^*)(f \w
1_\A))=$$$$=\text{\textbf{Sub}}_{\,\mathbb{B}(
\A,\,\A\,),\,\mathbb{B}(\A,f)}(\mathbb{B}(\A,f^*)(\mathbb{B}(\A,f)(1_\A))
=\text{\textbf{Sub}}_{\,\mathbb{B}(
\A,\,\A\,),\,\textbf{T}^\A_f}(\textbf{T}^\A_f(1_\A))\simeq$$$$\simeq
\text{\textbf{Sub}}_{\,\mathbb{B}(
\A,\,\A\,),\,{\textbf{T}^l}_{\mathbb{S}_f}}(\textbf{T}^l_{\mathbb{S}_f}(1_\A))\simeq
\textbf{I}^l_{\mathbb{B}(\A, \A)}(\mathbb{S}_f).
$$ Since the unit of the adjunction $$\mathbb{B}(\A,f) \dashv
\mathbb{B}(\A,f^*) : \mathbb{B}(\A, \A)\to \mathbb{B}(\A, \B)$$ is a
monomorphism, the morphism $\eta_f : 1 \to f^* \w f$ is right pure
in the monoidal category $\mathbb{B}(\A,\A)$, implying that
$\textbf{I}^l_{\mathbb{B}(\A, \A)}(\mathbb{S}_f)$ has the structure
of a monoid (as in Proposition \ref{monoide}). Transporting this
monoid structure along the above chain of bijections, one gets the
structure of a monoid on
$$\text{\textbf{Sub}}_{\,\mathbb{B}(
\A,\,\A\,),\,\mathbb{B}(\A,f)}(\mathbb{B}(\A,f^*)(f))$$ and it is
straightforward to check that this structure is just the one
described in the proposition.
\end{proof}

\bigskip

It follows from the proof of Proposition \ref{monoide2} that for any
adjunction $f \dashv f^* : \B \to \A$ in a 2-category $\mathbb{K}$,
$$\text{\textbf{Sub}}_{\,\mathbb{K}(
\A,\,\A\,),\,\mathbb{K}(\A,f)}(\mathbb{K}(\A,f^*)(f))=\textbf{I}^l_{\mathbb{K}(
\A,\,\A\,)}(S_f).$$

Consider now the composite
$$\Gamma^\mathbb{K}_f =\gamma' \Psi_{\mathbb{K}(\A,f),\,f} :\textbf{I}^l_{\mathbb{K}(
\A,\,\A\,)}(S_f)= \text{\textbf{Sub}}_{\mathbb{K}(\A,\A),\,
\mathbb{K}(\A,f)}(\mathbb{K}(\A,f^*)(f)) \to
\text{End}_\B(\mathfrak{C}_f,\mathfrak{C}_f).$$ This map takes $[(h,
i_h)]$ to the composite
$$\xymatrix{
ff^* \ar[rr]^-{(\xi_{i_h})^{-1}f^*}&&fhf^* \ar[rr]^-{f\eta_f
hf^*}&&ff^*fhf^* \ar[rr]^-{ff^*\xi_{i_h} f^* }&& ff^*ff^*
\ar[r]^-{ff^* \e_f}& ff^*}.$$ But since $\xi_{i_h}=\varepsilon_f f
\cdot f^*i_h$, it follows from the naturality of $\eta_f -$ that the
diagram
$$
\xymatrix{ h \ar[rr]^{\eta_f h} \ar[dd]_{i_h}&& f^*fh
\ar[rrdd]^{f^*\xi_{i_h} }\ar[dd]_{f^*fi_h}&&&&\\\\
f^*f \ar[rr]_{\eta_f f^*f}&& f^*ff^*f \ar[rr]_{f^* \varepsilon_f
f}&&f^*f}$$ commutes. Thus
$$f^*\xi_{i_h} \cdot \eta_f h=f^* \e_f f \cdot \eta_f f^*f \cdot i_h= i_h $$ and
hence
$$\Gamma^\mathbb{K}_f([(h, i_h)])= ff^* \varepsilon_f\cdot fi_h f^*\cdot (\xi_{i_h})^{-1}f^*.$$

By abuse of notation, we will write $\Gamma_f$ when $\mathbb{K}$
clear from the context.

We are now in position to state our main theorem in the framework of
$2$--categories. It will be extended later to bicategories.

\begin{theorem}\label{main2cat} When the functor $$\mathbb{K}(\A,f): \mathbb{K}(\A, \A) \to \mathbb{K}(\A, \B)$$
is comonadic, the map
\[
\Gamma^\mathbb{K}_f : \mathbf{I}^l_{\mathbb{K}( \A,\,\A\,)}(S_f) \to
\mathrm{End}_\B(\mathfrak{C}_f,\mathfrak{C}_f)
\]
 is an isomorphism of monoids.
\end{theorem}
\begin{proof} When the functor $$\mathbb{K}(\A,f) : \mathbb{K}(\A, \A) \to \mathbb{K}(\A,
\B)$$ is comonadic, it follows from Propositions
\ref{subobjetoscoalgebrasbiy} and \ref{endocoal} that the map
$\Gamma_f$ is bijective. So it only remains to show that $\Gamma_f$
is a morphism of monoids. We shall prove that the inverse
$\Gamma'_f$ to $\Gamma_f$ is a monoid isomorphism.

Note first that since the functor $\mathbb{K}(\A,f)$ is assumed to
be comonadic, the diagram
$$
\xymatrix{ 1_\A \ar[r]^-{\eta_f ^\A}& \mathbb{K}(\A,f^*f)
\ar@{->}@<0.5ex>[rr]^-{\eta_f
\mathbb{K}(\A,f^*f)}\ar@{->}@<-0.5ex>[rr]_-{\mathbb{K}(\A,f^*f)
\eta_f} && \mathbb{K}(\A,f^*ff^*f),}$$ and hence also its
$1_\A$-component
$$
\xymatrix{ 1 \ar[r]^-{\eta_f}&f^*f \ar@{->}@<0.5ex>[rr]^-{\eta_f
f^*f}\ar@{->}@<-0.5ex>[rr]_-{f^*f \eta_f} && f^*ff^*f}$$ is an
equalizer (see, for example, \cite{BW}). Thus
$\Gamma'_f(1_{\mathfrak{C}_f})=[\text{eq}(\eta_f f^*f,
f^*f\eta_f)]=[(1, \eta_f)]$, proving that $\Gamma'_f$ preserves the
identity element.

Next, according to \eqref{adj3}, the following diagram is
commutative:

\begin{equation}\label{adj5}
\xymatrix{ ff^*f \ar[dd]_{f \eta_f f^* f}\ar[rr]^-{f \eta_f f^*
f}\ar@{=}[ddrr]&&ff^*ff^*f
\ar[dd]^{ff^* \e_f f}\\\\
ff^*ff^*f \ar[rr]_{\e_f ff^*f}&& ff^*f,}
\end{equation}

Since $\Gamma'_f=\overline{\Psi}_{\mathbb{K}(\A,f),\,f}\cdot
\gamma,$ for any $\alpha \in
\text{End}_{\B}(\mathfrak{C}_f,\mathfrak{C}_f)$,
$$\Gamma'_f (\alpha)=\overline{\Psi}_{\mathbb{K}(\A,f),\,f}
(\gamma(\alpha))=\overline{\Psi}_{\mathbb{K}(\A,f),\,f}(f, \alpha
f\cdot f \eta_f)=[(f_\alpha, i_\alpha)],$$ where
$(f_\alpha,i_\alpha)$ is an equalizer of the pair

\begin{equation}\label{adj6}\xymatrix{ f^* f \ar@/_2pc/@{->}[rrrr]_{\eta_f f^*f }
\ar[rr]^{f^*f \eta_f} && f^* f f^* f \ar[rr]^{f^*  \alpha f} && f^*
f f^* f\,}\end{equation}

So, in particular,
\begin{equation}\label{adj7}
 \eta_f f^* f \cdot i_\alpha =f^* \alpha f\cdot (f^* f \eta_f)\cdot i_\alpha
\end{equation}

Suppose now that $\alpha_1, \alpha_2 \in
\text{End}_{\B}(\mathfrak{C}_f,\mathfrak{C}_f),$ and consider the
product
$$\Gamma'_f(\alpha_1) \cdot \Gamma'_f(\alpha_2)=[(f_{\alpha_1},
i_{\alpha_1})]\cdot[(f_{\alpha_2}, i_{\alpha_2})]=[(f_{\alpha_1}
f_{\alpha_2}, i_{f_{\alpha_1} f_{\alpha_2}})].$$

We have:
$$\eta_f f^* f \cdot f^* \e_f f \cdot i_{\alpha_1}i_{\alpha_2}=$$$$=\eta_f f^* f
\cdot f^* \e_f f \cdot i_{\alpha_1}f^*f \cdot h_1
i_{\alpha_2}=\,\,\text{by functoriality \,\,of \,\,composition}
$$$$=f^*ff^* \e_f f \cdot \eta_f f^*ff^*f \cdot i_{\alpha_1}f^*f \cdot h_1
i_{\alpha_2}=\,\,\text{by \,\,\eqref{adj7}}
$$$$=f^*ff^* \e_f f \cdot f^*\alpha_1 ff^*f \cdot f^*f \eta_f f^*f
\cdot i_{\alpha_1}f^*f \cdot h_1 i_{\alpha_2}=\,\,\text{by
functoriality \,\,of \,\,composition}
$$$$=f^* \alpha_1 f \cdot f^*ff^* \e_f f\cdot f^*f \eta_f f^*f
\cdot i_{\alpha_1}f^*f \cdot h_1 i_{\alpha_2}=\,\,\text{by
\,\,\eqref{adj5}}$$$$=f^* \alpha_1 f \cdot f^* \e_f ff^*  f\cdot
f^*f \eta_f f^*f \cdot i_{\alpha_1}f^*f \cdot h_1 i_{\alpha_2}=$$
$$=f^* \alpha_1 f \cdot f^* \e_f ff^*  f\cdot f^*f
\eta_f f^*f \cdot f^*f i_{\alpha_2} \cdot
i_{\alpha_1}h_2=\,\,\text{by \,\,\eqref{adj7}}$$$$=f^* \alpha_1 f
\cdot f^* \e_f ff^*  f \cdot f^*f f^* \alpha_2 f \cdot f^*f f f^*
\eta_f \cdot f^*f i_{\alpha_2} \cdot i_{\alpha_1}h_2=\,\,\text{by
functoriality \,\,of \,\,composition}
$$$$=f^* \alpha_1 f \cdot f^* \alpha_2 f \cdot f^*
\e_f ff^*  f \cdot f^*f f f^* \eta_f \cdot f^*f i_{\alpha_2} \cdot
i_{\alpha_1}h_2=$$$$=f^* (\alpha_1 \alpha_2) f \cdot f^* \e_f ff^* f
\cdot f^*f f f^* \eta_f \cdot f^*f i_{\alpha_2} \cdot
i_{\alpha_1}h_2=$$$$=f^* (\alpha_1 \alpha_2) f \cdot f^* \e_f ff^* f
\cdot f^*f f f^* \eta_f \cdot i_{\alpha_1}i_{\alpha_2}=\,\,\text{by
functoriality \,\,of \,\,composition}$$$$=f^* (\alpha_1 \alpha_2) f
\cdot f^*f\eta_f \cdot  f^* \e_f f \cdot i_{\alpha_1}i_{\alpha_2}.$$

Since $i_{f_{\alpha_1}f_{\alpha_2}}=f^* \e_f f \cdot
i_{\alpha_1}i_{\alpha_2} $, it follows that
$$f^*(\alpha_1 \alpha_2)f \cdot
f^*f \eta_f \cdot i_{f_{\alpha_1}f_{\alpha_2}}= \eta_f f^*f \cdot
i_{f_{\alpha_1}f_{\alpha_2}}
$$ and since $(f_{\alpha_1 \alpha_2}, i_{\alpha_1 \alpha_2})$ is
an equalizer of the diagram \eqref{adj6}, there exists a unique
morphism $k_{\alpha_1, \alpha_2}: f_{\alpha_1}f_{\alpha_2} \to
f_{\alpha_1 \alpha_2} $ making the diagram
$$
\xymatrix{
f_{\alpha_1}f_{\alpha_2} \ar[dd]_{k_{\alpha_1, \alpha_2}} \ar[rrd]^{i_{f_{\alpha_1} f_{\alpha_2}}}&&\\
&& f^* f\\
f_{\alpha_1 \alpha_2} \ar[rru]_{i_{\alpha_1 \alpha_2}}&&}$$ commute.
It now follows -since any comonadic functor is conservative- from
Lemma \ref{conservativosubobjetos} that $[i_{f_{\alpha_1}
f_{\alpha_2}}]=[i_{\alpha_1 \alpha_2}]$ in
$\text{\textbf{Sub}}_{\mathbb{K}(\A,\A),\,\mathbb{K}(\A,f)}(f^* f).$
Therefore
$$\Gamma'_f (\alpha_1 \alpha_2)=\Gamma'_f (\alpha_1)\cdot
\Gamma'_f (\alpha_2).
$$ It clearly implies that $\Gamma_f$ is an isomorphism of monoids. This completes the proof.
\end{proof}

Next, we give the dual statement to Theorem \ref{main2cat}.

\begin{theorem}\label{main2catbis} Let $\eta_f, \e_f :
f \dashv f^* : \B \to \A$ be an adjunction in a 2-category
$\mathbb{K}$ such that the functor $$\mathbb{K}(f^*,\A):
\mathbb{K}(\A, \A) \to \mathbb{K}(\B, \A)$$ is comonadic. Then the
assignment
$$\xymatrix{[(h, i_h)]\longrightarrow
(ff^* \ar[rr]^-{f(\xi_{i_h})^{-1}}&&fhf^* \ar[rr]^-{fi_hf^* }&&
ff^*ff^* \ar[r]^-{ff^* \e_f}& ff^*)}$$ yields an isomorphism of
monoids
$$\Gamma^\mathbb{K}_{f^*}=\Gamma_{f^*}:\mathbf{I}^r_{\mathbb{K}(\A, \A)}(S_f) \to
(\text{\emph{End}}_\B(\mathfrak{C}_f,\mathfrak{C}_f))^{\text{op}}.$$
\end{theorem}

\begin{proof}Since $f^*$ is right adjoint to $f$ in $\mathbb{K}$,
$f$ is right adjoint to $f^*$ in $\mathbb{K}^t$. Thus the adjunction
$\mathbb{K}(f^*,\A)\dashv \mathbb{K}(f,\A) : \mathbb{K}(\B, \A) \to
\mathbb{K}(\A, \A)$ in $\mathbb{K}$ can be seen as the adjunction
$\mathbb{K}^t(\A,f^*)\dashv \mathbb{K}^t(\A,f) : \mathbb{K}^t(\A,
\B) \to \mathbb{K}^t(\A, \A)$ in $\mathbb{K}^t$. And since the
functor $\mathbb{K}(f^*,\A)$ is comonadic by assumption, the functor
$\mathbb{K}^t(\A,f^*)$ is also comonadic, and thus we can apply
Theorem \ref{main2cat} to the adjunction $\mathbb{K}^t(\A,f^*)\dashv
\mathbb{K}^t(\A,f)$ to conclude that the map
$$\Gamma^{\mathbb{K}^t}_{f^*}: \mathbf{I}^l_{\mathbb{K}^t(\A, \A)}
(S^{\mathbb{K}^t}_{f^*})\to
\text{End}_\B(\mathfrak{C}^{\mathbb{K}^t}_{f^*},\mathfrak{C}^{\mathbb{K}^t}_{f^*})$$
is an isomorphism of monoids. Now, using the fact that
$\mathbb{K}^t(\A, \A)$ can be identified with the monoidal category
$(\mathbb{K}(\A, \A))^t$, it is not hard to show that
$$\mathbf{I}^l_{\mathbb{K}^t(\A, \A)}
(S^{\mathbb{K}^t}_{f^*})=\mathbf{I}^r_{\mathbb{K}(\A, \A)}
(S^{\mathbb{K}}_{f}),
\,\,\,\text{End}_\B(\mathfrak{C}^{\mathbb{K}^t}_{f^*},\mathfrak{C}^{\mathbb{K}^t}_{f^*})
=(\text{End}_\B(\mathfrak{C}_f,\mathfrak{C}_f))^{\text{op}},$$ and
that $\Gamma^{\mathbb{K}^t}_{f^*}$ is just the map
$\Gamma^\mathbb{K}_{f^*}$. Consequently, $\Gamma^\mathbb{K}_{f^*}$
is an isomorphism of monoids.
\end{proof}

Our next aim is to lift Theorem \ref{main2cat} from $2$--categories
to bicategories. We will proceed step by step.

\begin{proposition}\label{reduccionaestricto2} Let $\mathbb{B}$ be a bicategory, $\mathbb{K}$ be a 2-category, and
$\phi=(\phi, \overline{\phi}, \phi_0):\mathbb{B} \to \mathbb{K}$ be
a homomorphism of bicategories. Then for any adjunction $\eta_f,
\e_f : f \dashv f^* : \B \to \A$ in $\mathbb{K}$, the (iso)morphism
$\overline{\phi}_{f^*, f}: \phi(f^*) \phi(f) \to \phi (f^* \w f)$
yields an isomorphism of $\mathbb{K}(\phi(\A), \phi(\A))$-monoids
$$\overline{\kappa}_f^\phi: S^\mathbb{K}_{\phi(f)} \to
\phi(S^\mathbb{B}_f).$$
\end{proposition}

\begin{proof} Considering the diagram
$$\xymatrix@R=19pt@C=4pt{\phi(f^*)\phi(f)\phi(f^*)\phi(f)
\ar[rrrr]^-{\overline{\phi}\phi(f^*)\phi(f)}
\ar[dddd]_{\phi(f^*)\overline{\phi}\phi(f)}&&&& \phi(f^*\w f)
\phi(f^*)\phi(f) \ar[dd]^{\overline{\phi}
\phi(f)}\ar[rrrr]^{\phi(f^*\w f)\overline{\phi}}&&&&
\phi(f^*\w f)\phi(f^*\w f)  \ar[dd]^{\overline{\phi}}\\\\
&&&& \phi((f^* \w f)\w f^*)\phi(f)\ar@{}[rr]^-{(2)}
\ar[dd]_{\phi(\alpha_{f^*, f, f^*})\phi(f)}
\ar[rrrrdd]_{\overline{\phi}}&&&& \phi((f^* \w f)\w (f^* \w f)) \\\\
\phi(f^*)\phi(f \w f^*)\phi(f)\ar@{}[rruuuu]_{(1)}
\ar[rrrr]^{\overline{\phi}\phi(f)
}\ar[dd]_{\phi(f^*)\phi(\varepsilon_f)\phi(f)}&&&& \phi(f^* \w (f \w
f^*))\phi(f) \ar@{}[rr]^-{(3)} \ar[dd]_{\phi(f^* \w
\varepsilon_f)\phi(f)} \ar[rrrrdd]_{\overline{\phi}}&&&&
\phi(((f^* \w f)\w f^*)\w f) \ar[dd]^{\phi(\alpha_{f^*, f, f}\w f)} \ar[uu]_{\phi(\alpha_{f^* \w f, f^*, f})}\\\\
\phi(f^*)\phi(1_\A)\phi(f) \ar@{}[rruu]_-{(4)}
\ar[dd]_{\phi(f^*)(\phi_0)^{-1}_\A\phi(f)}\ar[rrrr]_{\overline{\phi}\phi(f)}&&&&
\phi(f^* \w 1_\A)\phi(f) \ar@{}[rr]^-{(5)}
\ar[dd]_{\phi(r_{f^*})\phi(f)} \ar[rrrrdd]_{\overline{\phi}}&&&&
\phi((f^* \w( f\w f^*))\w f) \ar[dd]^{\phi((f^* \w \varepsilon_f)\w f)}\\\\
\phi(f^*)1_{\phi(\A)}\phi(f) \ar@{=}[rrrr]\ar@{}[rruu]_-{(6)}&&&&
\phi(f^*)\phi(f)\ar[rrrrdd]_{\overline{\phi}}\ar@{}[rr]^-{(7)}&&&&
\phi((f^* \w 1_\A)\w f)\ar[dd]^{\phi(r_{f^*}\w f)}\\\\
&&&&&&&& \phi(f^* \w f)}
$$ in which

\begin{itemize}
\item Diagrams (1), (2) and (6) commute because $\phi$ is a
homomorphism of bicategories;

\item Diagrams (3), (4), (5) and (7) commute by naturality of
$\overline{\phi}$,
\end{itemize} we see that $$m_{\phi(S^{\mathbb{B}}_f)}\cdot (\overline{\phi}_{f^*, f}\overline{\phi}_{f^*,
f})=$$$$=(\phi((\phi_0)_\A)\w f)\cdot (\phi((f^* \w \varepsilon_f)\w
f))\cdot (\phi(\alpha_{f^*, f, f}\w f))\cdot (\phi(\alpha^{-1}_{f^*
\w f, f^*, f}))\cdot (\overline{\phi}_{f^* \w f, f^*, f})\cdot
(\overline{\phi}_{f^*, f}\overline{\phi}_{f^*,
f})=$$$$=\overline{\phi}_{f^*, f} \cdot
(\phi(f^*)(\phi_0)^{-1}_\A\phi(f))\cdot
(\phi(f^*)\phi(\varepsilon_f)\phi(f))\cdot
(\phi(f^*)\overline{\phi}_{f, f^*}\phi(f))=\overline{\phi}_{f^*, f}
\cdot m_{S^{\mathbb{K}}_{\phi(f)}}.$$ Moreover, since
$e_{\phi(S^{\mathbb{B}}_f)}=\phi(\eta_f)\cdot (\phi_0)_\A
:1_{\phi(\A)}\to \phi(f^* \w f)$, the diagram
$$\xymatrix{
1_{\phi(\A)}\ar[d]_{(\phi_0)_\A}
\ar[rrrdd]^{e_{\phi(S^{\mathbb{B}}_f)}}&&&\\
\phi(1_\A) \ar[d]_{\phi(\eta_f)}&&&\\
\phi(f^* \w f) \ar[rr]_-{(\overline{\phi}_{f^*, f})^{-1}}&&
\phi(f^*)\phi(f) \ar[r]_{\overline{\phi}_{f^*, f}}& \phi(f^* \w f)}
$$ clearly commutes, showing that $\overline{\phi}_{f^*, f}$
preserves the unit. Thus $\overline{\phi}_{f^*, f}: \phi(f^*)
\phi(f) \to \phi (f^* \w f)$ is an isomorphism of monoids in
$\mathbb{K}(\phi(\A), \phi(\A)).$
\end{proof}

Dually, one has:

\begin{proposition} Let $\mathbb{B}$ be a bicategory,
 $\mathbb{K}$ be a 2-category, and
$\phi=(\phi, \overline{\phi}, \phi_0):\mathbb{B} \to \mathbb{K}$ be
a homomorphism of bicategories. Then for any adjunction $\eta_f,
\e_f : f \dashv f^* : \B \to \A$ in $\mathbb{K}$, the (iso)morphism
$\overline{\phi}_{f^*, f}: \phi(f^*) \phi(f) \to \phi (f^* \w f)$
yields an isomorphism of $\mathbb{K}(\phi(\A), \phi(\A))$-comonoids
$$\underline{\kappa}_f^\phi: \mathfrak{C}^\mathbb{K}_{\phi(f)} \to
\phi(\mathfrak{C}^\mathbb{B}_f).$$
\end{proposition}

We are now ready to prove our main result, namely, the following
bicategorical version of Masuoka's Theorem.

\begin{theorem}\label{maincomonadic} Let $\eta_f, \e_f : f \dashv f^*  : \B \to \A$
be an adjunction in a bicategory $\mathbb{B}$ such that the functor
$\mathbb{B}(\A, f):\mathbb{B}(\A, \A) \to \mathbb{B}(\A, \B)$ is
comonadic. Then the rule $$[(h, i_h)]\longrightarrow \theta_h,$$
where $\theta_h$ is the composite
$$(f \w r_{f^*})
\cdot (f\w (f^* \w \varepsilon_f)) \cdot (f \w \alpha_{f^*, f, f^*})
\cdot (\alpha_{f, f^* \w f, f^*})\cdot ((f \w i_h)\w
f^*)\cdot((\xi_{i_h})^{-1}\w f^*),$$  defines an isomorphism of
monoids:
$$\Gamma^\mathbb{B}_{f}=\Gamma_{f}:\mathbf{I}^l_{\mathbb{B}(\A,
\A)}(S^{\mathbb{B}}_f) \to
\text{\emph{End}}_\B(\mathfrak{C}^{\mathbb{B}}_f,\mathfrak{C}^{\mathbb{B}}_f).$$
\end{theorem}
\begin{proof}According to \cite{L}, there exist a 2-category
$\mathbb{K}$, homomorphisms $$\phi=(\phi, \overline{\phi}, \phi_0):
\mathbb{B} \to \mathbb{K}, \,\,\phi'=(\phi', \overline{\phi'},
\phi'_0): \mathbb{K} \to \mathbb{B}$$ of bicategories and invertible
icons $\phi\phi'\simeq 1$ and $\phi' \phi \simeq 1$. And considering
the monoidal functor
$$\phi_{\A, \A}:\mathbb{B}(\A, \A) \to \mathbb{K}(\phi(\A),
\phi(\A)),$$ one observes that the composite
$$\xymatrix{
\mathbf{I}^l_{\mathbb{B}(\A, \A)}(S^{\mathbb{B}}_f)
\ar[rr]^-{(\phi_{\A, \A})_{S^{\mathbb{B}}_f}}&&
\mathbf{I}^l_{\mathbb{K}(\phi(\A), \phi(\A))}(\phi_{\A,
\A}(S^{\mathbb{B}}_f))\ar[rrr]^-{\mathbf{I}^l_{\mathbb{K}(\phi(\A),
\phi(\A))}((\overline{\kappa}_f^\phi)^{-1})} &&&
\mathbf{I}^l_{\mathbb{K}(\phi(\A),
\phi(\A))}(S^{\mathbb{K}}_{\phi(f)})}$$
$$\xymatrix{\ar[r]^-{\Gamma^{\mathbb{K}}_{\phi(f)}}&
\text{{End}}_{\phi(\B)}(\mathfrak{C}^{\mathbb{K}}_{\phi(f)},
\mathfrak{C}^{\mathbb{K}}_{\phi(f)})\ar[rr]^-{t^{\underline{\kappa}^\phi_f}_{\mathbb{K}}}&&
\text{{End}}_{\phi(\B)}(\phi(\mathfrak{C}^{\mathbb{B}}_{f}),
\phi(\mathfrak{C}^{\mathbb{B}}_{f}))\ar[rr]^-{(s^\phi_{\mathfrak{C}^{\mathbb{B}_f}})^{-1}}&&
\text{{End}}_{\B}(\mathfrak{C}^{\mathbb{B}}_{f},\mathfrak{C}^{\mathbb{B}}_{f})}$$
is just the map $\Gamma^\mathbb{B}_{f}$. Now, since $\phi_{\A, \A}$
is a strong monoidal functor, it follows that
\begin{itemize}

\item $(\phi_{\A, \A})_{S^{\mathbb{B}}_f}$ is an isomorphism of
monoids by  Proposition \ref{reduccionaestricto};

\item $\mathbf{I}^l_{\mathbb{K}(\phi(\A),
\phi(\A))}(((\overline{\kappa}_f^\phi)^{-1}))$ is an isomorphism of
monoids by Proposition \ref{reduccionaestricto2};

\item $t^{\underline{\kappa}^\phi_f}_{\mathbb{K}}$ is an
isomorphism of monoids by Proposition \ref{isoiso}(1);

\item $s^\phi_{\mathfrak{C}^{\mathbb{B}_f}}$ is an isomorphism of
monoids by Proposition \ref{isoiso} (2).
\end{itemize} Moreover, considering the following diagram
$$\xymatrix{
\mathbb{B}(\A,\A) \ar[d]_{\phi_{\A,\,\A}}
\ar@{->}@<0.5ex>[rr]^-{\mathbb{B}(\A, \,f)} &&
\mathbb{B}(\B,\A) \ar@{->}@<0.5ex>[ll]^-{\mathbb{B}(\A, \,f^*)} \ar[d]^{\phi_{\B,\,\A}}\\
\mathbb{K}(\phi(\A),\phi(\A))
\ar@{->}@<0.5ex>[rr]^-{\mathbb{K}(\phi(\A),\,\phi(f))} &&
\ar@{->}@<0.5ex>[ll]^-{\mathbb{K}(\phi(\A),\,\phi(f^*))}
\mathbb{B}(\phi(\A),t(\B))}$$ which is commutative (up to
isomorphism), since $\phi$ is a homomorphism of bicategories, and
using the fact that $\phi_{\A, \A}$ and $\phi_{\A, \B}$ are both
equivalences of categories, one sees that the functor
$\mathbb{K}(\phi(\A), \phi(f)): \mathbb{K}(\phi(\A), \phi(\A)) \to
\mathbb{K}(\phi(\A), \phi(\B))$ is comonadic, and applying Theorem
\ref{main2cat} gives that $\Gamma^{\mathbb{K}}_{\phi(f)}$ is an
isomorphism of monoids. Consequently, the map
$\Gamma^\mathbb{B}_{f}$ is an isomorphism of monoids.
\end{proof}

Dually, one has:

\begin{theorem}\label{maincomonadicdual}
Let $\eta_f, \e_f : f \dashv f^* : \B \to \A$ be an adjunction in a
bicategory $\mathbb{B}$ such that the functor $\mathbb{B}( f^*,\,
\A):\mathbb{B}(\A, \A) \to \mathbb{B}(\B, \A)$ is comonadic. Then
there is an isomorphism of monoids:
$$\Gamma^\mathbb{B}_{f^*}=\Gamma_{f}:\mathbf{I}^r_{\mathbb{B}(\A,
\A)}(S^{\mathbb{B}}_f) \to
(\text{\emph{End}}_\A(\mathfrak{C}^{\mathbb{B}}_f,\mathfrak{C}^{\mathbb{B}}_f))^{op}.$$
\end{theorem}

\section{Sufficient conditions for the comonadicity: flatness, purity and
separability}\label{flatpurityseparability}

We continue to suppose that $\mathbb{B}$ is a bicategory. In this
section, we find sufficient conditions on the $1$--cell $f$ that
imply the comonadicity of the functor $\mathbb{B}(\mathcal{C},f)$,
and thus the existence of the isomorphism of monoids $\Gamma_f$.
 We start with
a definition inspired by the (bi)module case.

\begin{definition}Let $\C$ be an arbitrary 0-cell of $\mathbb{B}$.
\begin{itemize}
\item [\emph{(i)}] A 1-cell $f: \A \to \B$ is right (resp. left)
$\C$-\emph{flat} if the functor $$\mathbb{B}(\C, f)=f\w -
:\mathbb{B}(\C, \A)\to \mathbb{B} (\C, \B) $$ (resp.
$$\mathbb{B}(f, \C)= -\w f :\mathbb{B}(\B, \C)\to \mathbb{B}(\A,
\C)) $$ preserves equalizers. \item [\emph{(ii)}] A 2-cell $\tau :f
\to g : \A \to \B$ is \emph{right} ( resp. \emph{left})
$\C$-\emph{pure} if for any 1-cell $h: \C \to \A$ (resp. $h: \B \to
\C$) the morphism $\tau \w 1_h : f \w h \to g \w h$ (resp. $1_h \w
\tau : h\w f \to h \w g$ ) is a regular monomorphism in the category
$\mathbb{B}(\C, \B)$ (resp. $\mathbb{B}(\A, \C)$).
\end{itemize}
\end{definition}

A consequence of \cite[Theorem 3.10]{EG1} is that if $M$ is an
$A$--$B$--bimodule over rings $A,B$ such that $M_B$ is finitely
generated and projective and ${}_AM$ is faithfully flat, then $-
\otimes_A M : \mathrm{Mod}_A \rightarrow \mathrm{Mod}_B$ is
comonadic. This statement could be alternatively deduced from Beck's
comonadicity theorem. The following proposition makes it clear in
the framework of our general bicategory $\mathbb{B}$.

\begin{proposition}\label{flatpure} Let $\eta_f,\e_f : f \dashv f^* : \B \to \A$ be
an adjunction in $\mathbb{B}$. Suppose that for a 0-cell $\C$, $f$
(resp. $f^*$) is right $\C$-flat (resp. $f^*$ is left $\C$-flat) and
$\eta_f$ is right (resp. left) $\C$-pure. Then the functor
$$\mathbb{B}(\C, f)=f \w - :\mathbb{B}(\C, \A) \to \mathbb{B}(\C, \B) $$ (resp. $$
\mathbb{B}(f^*, \C)=-\w f^*  :\mathbb{B}(\A, \C)\to \mathbb{B}(\B,
\C)) $$ is comonadic.
\end{proposition}
\begin{proof} By duality, it suffices to prove the first statement.
The functor $$\mathbb{B}(\C, f^*)=f^* \w - :\mathbb{B}(\C, \B) \to
\mathbb{B}(\C, \A)$$ is right adjoint to the functor
$$\mathbb{B}(\C, f)=f \w -:\mathbb{B}(\C, \A) \to \mathbb{B}(\C,
\B) $$ and for any 1-cell $h: \C \to \A,$ the $h$-component of the
unit $\eta_{\C}^f$ of this adjunction is given by the composite
$$\alpha_{f^*, \,f, \,f^*}\cdot(\eta_f \w h)\cdot (l_h)^{-1}: h \to f^* \w f \w h.$$ Since
$\eta_f$ is right $\C$-pure by hypothesis, $\eta_f \w h$ (and hence
also the composite $\alpha_{f^*, \,f, \,f^*}\cdot(\eta_f \w h)\cdot
(l_h)^{-1}$) is a regular monomorphism. Hence the unit of the
adjunction $\mathbb{B}(\C, f) \dashv \mathbb{B}(\C, f^*)$ is a
regular monomorphism, implying that the functor $\mathbb{B}(\C, f)$
reflects isomorphisms (see, for example, \cite{A}). Moreover, since
$f$ is right $\C$-pure, the functor $\mathbb{B}(\C, f)$ preserves
equalizers, and hence the functor $\mathbb{B}(\C, f)$ is comonadic
by a simple and well-known application of (the dual of) Beck's
theorem (see, \cite{M}).
\end{proof}

We say that a 1-cell $f : \A \to \B$ is right (resp. left) flat if
for any 0-cell $\C$, $f$ is right (resp. left) $\C-$flat. Similarly,
we say that a 2-cell $\tau : f \to f': \A \to \B$ is right (resp.
left) pure if for any 0-cell $\C$, $f$ is right (resp. left)
$\C-$pure.

The following is an immediate consequence of Proposition
\ref{flatpure}.

\begin{proposition}\label{planoimplicacomonadico}
Let $\eta_f,\e_f : f \dashv f^* : \B \to \A$ be an adjunction in
$\mathbb{B}$ and suppose that $f$ is right (resp. $ f^*$ is left)
flat and that $\eta_f$ is right (resp. left) pure. Then, for any
0-cell $\C$, the functor
$$\mathbb{B}(\C, f) :\mathbb{B}(\C, \A) \to \mathbb{B}(\C, \B) $$
(resp.$$\mathbb{B}(f^*, \C) :\mathbb{B}(\A, \C)\to \mathbb{B}(\B,\C)
)$$ is comonadic.
\end{proposition}

In the light of Proposition \ref{planoimplicacomonadico}, we get
from Theorems \ref{main2cat} and \ref{main2catbis} the following
generalization of Masuoka's theorem.

\begin{theorem}\label{flatpureiso} Let $\eta_f,\e_f : f \dashv f^* : \B \to \A$
be an adjunction in $\mathbb{B}$.
\begin{itemize}
\item [\emph{(i)}] If $f$ is right $\A$-flat and $\eta_f$ is right
$\A$-pure, then the map $$\Gamma_f :
{\textbf{\emph{I}}}^l_{\mathbb{B}(\A,\,\A)}(\mathbb{S}_f) \to
\text{\emph{End}}_\B(\mathfrak{C}_f,\mathfrak{C}_f)$$ is an
isomorphism of monoids.

\item [\emph{(ii)}] If $f$ is left $\A$-flat and $\eta_f$ is left
$\A$-pure, then the map
$$\Gamma_{f^*}:({\textbf{\emph{I}}}^r_{\mathbb{B}(\A,\,\A)}(\mathbb{S}_f))^{\text{op}}\to
\text{\emph{End}}_\B(\mathfrak{C}_f, \mathfrak{C}_f)$$ is an
isomorphism of monoids.
\end{itemize}
\end{theorem}

We shall need the following

\begin{proposition}\label{purepures}
 Let $\V=(V, \otimes, I, \alpha, l,
r)$ be a monoidal category and consider a monoid $\mathbb{S} = (S,
e_S, m_S)$ over $\V$. If the morphism $e_S : I \to S$ is right pure,
then for any $[(J, i_J)]\in \textbf{I}^l_\V(S)$, the morphism $i_J :
J \to S$ is also right pure.
\end{proposition}
\begin{proof} Since the morphism $e_S$ is right pure, $e_S \otimes X : I\w X \to S \w X$,
and hence also $$(e_S \otimes X)\cdot (l_{X})^{-1} :  X \to S \w
X,$$ is a monomorphim for any object $X$. Then in particular, the
morphism
$$(e_S \w J \w X)\cdot (l_{J \w X})^{-1} : J \w X \to S
\w J\w X $$ is a monomorphism, and since $[(J, i_J)]\in
\textbf{I}^l_\V(S)$ and thus $\xi_{i_J}$ is an isomorphism, the
composite
$$
\xymatrix{ J \w X \ar[rr]^-{(l_{J \w X})^{-1}}&& I \w J \w X
\ar[rr]^-{e_S \w J \w X }&& S \w J\w X \ar[rr]^-{\xi_{i_J} \w X}&& S
\w X }$$ is also a monomorphism. We claim that this composite is
just the morphism $i_J \w X$. Indeed, consider the following diagram
$$
\xymatrix{ J \w X \ar@{}[rrdd]^{(1)}\ar[rr]^-{(l_J)^{-1}\w X}\ar[dd]
_{i_J \w X}&& I \w J \w X \ar@{}[rrdd]^-{(2)} \ar[rr]^-{e_S \w J \w
S} \ar[dd]^{I \w i_J \w X}
&& S \w J \w X  \ar[dd]^-{S \w i_J \w  X} \\\\
S \w X \ar[rr]_-{(l_{S})^{-1} \w X} && I \w S \w X
\ar[rrdd]_{l_{S}\w X}\ar[rr]_-{e_S \w S \w X}
&& S \w S \w X \ar[dd]^{m_S \w X}\\\\
&&&& S \w X\,.}$$ In this diagram
\begin{itemize}
  \item Diagram (1) commutes by naturality of $l$;
  \item Diagram (2) commutes by bifunctoriality of $-\w -$, and
  \item the triangle commutes since $e_S: I \to S$ is
  a unit for the monoid $\mathbb{S}$.
\end{itemize} It follows $-$since $\xi_{i_J} = m_S\cdot (S \w i_J)$$-$
that
$$(m_S \w X)\cdot (S \w i_J \w  X)\cdot (e_S \w J \w X )
\cdot (l_{J \w X})^{-1}= \,\,\,\,\text{since\,\,}(l_{J \w X})^{-1}=
(l_J)^{-1}\w X) $$
$$=(m_S \w X)\cdot (S_S \w i_J \w  X)\cdot (e_S \w J \w X )
\cdot ((l_J)^{-1}\w X)=i_J \w X.$$ Therefore $i_J \w X$ is a
monomorphism for all $X \in \V $, proving that $i_J$ is right pure.
\end{proof}

\bigskip

An $\mathcal{A}$--ring $\mathbb{S} = (S,e_S,m_S)$ is said to be
\emph{split} if $e_S$ is a split monomorphism in
$\mathbb{B}(\mathcal{A},\mathcal{A})$. We will prove that an
$\mathcal{A}$--ring that comes from an adjunction is split if and
only if the representatives are separable functors in the sense of
\cite{NVV} (see \cite{rafael} for a characterization of separable
adjoint functors). Given an adjunction $\eta_f,\e_f : f \dashv f^* :
\B \to \A$ in a bicategory $\mathbb{B}$, one says that $f$ is a
\emph{separable 1-cell} (or that the 1-cell $f$ is \emph{separable})
if the 2-cell $\eta_f :1_\A \to f^* \w f$ is a split monomorphism in
the category $\mathbb{B}(\A, \A)$.

\begin{proposition}\label{separablesplit}
For an adjunction $\eta_f,\e_f : f \dashv f^* : \B \to \A$ in
$\mathbb{B}$, the following are equivalent:
\begin{itemize}
  \item [\emph{(i)}]$f$ is separable.
  \item [\emph{(ii)}] the $\A$-ring $S_f$ is split;
  \item [\emph{(iii)}] For each 0-cell $\C$, the unit $\eta^{\C}_f$ of
  the adjunction $$ \mathbb{B}(\C, f)\dashv \mathbb{B}(\C, f^*) :\mathbb{B}(\C, \B) \to \mathbb{B}
  (\C, \A) $$ is a split monomorphism (i.e., $\mathbb{B}(\C, f)$ is separable) ;
  \item [\emph{(iv)}]For each 0-cell $\C$, the unit $\overline{\eta}^{\C}_f$
  of the adjunction $$\mathbb{B}(f^*, \C) \dashv \mathbb{B}(f, \C) :\mathbb{B}(\B, \C) \to
  \mathbb{B}(\A, \C) $$ is a split monomorphism (i.e. $\mathbb{B}(f^*, \C)$ is separable).
\end{itemize}
\end{proposition}
\begin{proof} Since the unit of the $\A$-ring $S_f$ is just $\eta_f$, (i)
and (ii) are equivalent. Next, since
$(\eta^\C_{f})_f=\alpha_{f^*,\,f,\, f^*}\cdot(\eta_f \w f)\cdot
(l_f)^{-1}$ for all 1-cell $f: \C \to \A$, that (i) and (iii) are
equivalent is easily seen, while putting $\C=\A$ and $f=1_\A$ in
$(\eta^\C_{f})_f$ gives the implication (iii) $\Rightarrow$ (i).
Similarly, one can prove that (i) is equivalent to (iv).
\end{proof}

\begin{proposition}\label{separable} If $\eta_f,\e_f : f \dashv f^* : \B \to \A$ in an
adjunction $\mathbb{B}$ such that $f$ is a separable 1-cell and if
idempotents split in the category $\mathbb{B}(\A, \A)$, then the
functors $$\mathbb{B}(\C, f)=f \w - :\mathbb{B}(\C, \A) \to
\mathbb{B}(\C, \B)
$$ and $$\mathbb{B}(f^*, \C)=- \w f^*:\mathbb{B}(\A, \C) \to \mathbb{B}(\B, \C) $$
are both comonadic.
\end{proposition}
\begin{proof}Since $f$ is a separable 1-cell in $\mathbb{B}$ by hypothesis, it
follows from Proposition \ref{separablesplit} that the units of the
adjunctions
$$f \w - \dashv f^* \w - :\mathbb{B}(\C, \B) \to \mathbb{B}(\C,
\A)
$$ and $$- \w f^* \dashv - \w f :\mathbb{B}(\B, \C) \to
\mathbb{B}(\A, \C) $$ are both  split monomorphisms. And since
idempotents split in in the category $\mathbb{B}(\A, \A)$, one may
apply Proposition 3.16 of \cite{Me} to conclude that the functors
$\mathbb{B}(\C, f)$ and $\mathbb{B}(f^*, \C)$ are both comonadic.
\end{proof}

Combining  Proposition \ref{separable} and Theorems
\ref{maincomonadic} and \ref{maincomonadicdual}, we get:

\begin{theorem}\label{separableiso} Let  $\eta_f,\e_f :
f \dashv f^* : \B \to \A$ be an adjunction in a bicategory
$\mathbb{B}$ with $f$ a separable 1-cell. If idempotents split in
the category $\mathbb{B}(\A,\A)$, then each of the maps
$$\Gamma_f : \textbf{\emph{I}}^l_{\mathbb{B}(\A,\, \A)}(S_f) \to
\text{\emph{End}}_\B(\mathfrak{C}_f,\mathfrak{C}_f)$$ and
$$\Gamma_{f^*}:(\textbf{\emph{I}}^r_{\mathbb{B}(\A,\, \A)}(S_f))^{\text{op}}\to
\text{\emph{End}}_\B(\mathfrak{C}_f, \mathfrak{C}_f)$$ is an
isomorphism of monoids.
\end{theorem}

\section{Applications}\label{applications}

We will now consider bicategories of generalized bimodules in the
frameworks of functor additive categories and of firm modules, and
we will apply the theory so far developed to these cases. In
particular, the generalizations of Masuoka's theorem given in
\cite{EG} and \cite{Me1} will be deduced.

 Let $K$ be a
commutative ring with unit, and $\R$ the category of unital
$K$-modules. A $K$-category is a category $\A$ equipped with a
$K$-module structure on each hom set in such a way that composition
induces $K$-module homomorphisms $$\A (a, b) \otimes_K \A (b, c) \to
\A (a, c).$$

Equivalently, a $K$-category is a category which is enriched over
the symmetric monoidal closed category $\R$. A one-object
$K$-category is just an associative $K$-algebra. A $\mathbb
Z$-category is a preadditive category.

If $\A$ and $\B$ are $K$-categories, a functor $F : \A \to \B$ is
called a $K$-functor if for all $a, a' \in \A$, the map $$F_{a, a'}
: \A (a, a') \to \B (F(a), F(a'))$$ is a morphism of $K$-modules. A
$\mathbb Z$-functor is an additive functor. And when $\A$ is small,
we write $[\A, \B]$ for the category of $K$-functors from $\A$ to
$\B$, where $[\A, \B](F, F')$ is the $K$-module of $K$-natural
transformations from $F$ to $F'$.

Given $F : \A \to \R,$ $a, a' \in \A$, $f \in \A(a, a')$ and $x \in
\A(a),$ write $xf$ for $F(f)(x)$. Dually, if $G: \A^{\text{op}}\to
\R$ and $y \in \A(a'),$ write $fy$ for $F(f)(y)$.

If $\A$ and $\A'$ are $K$-categories, we let $\A \otimes_K \A'$
denote the $K$-category whose class of objects is $Ob(\A) \times
Ob(\A')$, where the $K$-module of morphisms from $(a, a')$ to $(b,
b')$ is $\A(a, b)\otimes_K \A(a', b')$.

Let $\A$ be a small $K$-category. It is well known (see \cite{F})
that there is a functor
$$- \otimes - : [\A, \R] \times [\A ^{op}, \R] \to
\R$$
$$(F, G) \to F \otimes G,$$
called the tensor product which is defined by the isomorphism
$$\R(F \otimes G, a)\simeq [\A ^{op}, \R](G, \R(F, a))$$
natural in $F \in [\A, \R]$, $G \in [\A ^{op}, \R]$, and $a \in \R$.
Here $\R(H, a)$ denotes the $K$-functor $\A ^{op} \to \R$ whose
value at $a'$ is $\R(H(a'), b)$. Recall that the tensor product $F
\otimes G$ is defined by $$F \otimes G=\bigoplus_{a \in \A} F(a)
\otimes_K G(a) /Q\,,$$ where $Q$ is the $K$-submodule of
$\bigoplus_{a \in \A} F(a) \otimes_K G(a) $ generated by all
elements of the form $(af \otimes_K a'-a \otimes_K fa')$.

Note that, if $F \in [\A ^{op}, \R]$ and $G \in [\A, \R]$, then $F
\otimes G \simeq G \otimes F$.

If $\A$ and $\B$ are small $K$-categories, an $\A$-$ \B$-bimodule
$\phi$ (which we shall denote by $ {\phi}: \A \rightsquigarrow \B $)
is a $K$-functor $\phi : \B ^{op} \otimes \A \to \R.$ Morphisms
between bimodules are just $K$-natural transformations. We shall
write $\text{\textbf{Bim}}(\A,\B)$ for the category of $\A$ -
$\B$-bimodules. If $\K$ is the unit $K$-category (that is, if $\K$
is the $K$-category with one object $*$ and with $\K(*, *)=K)$, a
bimodule $  \K \rightsquigarrow \A$ is essentially a $K$-functor $\A
^{op} \to \R$, and $ \text{\textbf{Bim}}(\K, \A)$ is (isomorphic to)
the category $[\A ^{op}, \R]$. Similarly, a bimodule $ \A
\rightsquigarrow \K $ is a $K$-functor $\A \to \R$ and
$\text{\textbf{Bim}}(\A,\K)$ is (isomorphic to) the category $[\A,
\R]$. And a bimodule $\K \rightsquigarrow \K $ can be identified
with a $K$-module and $\text{\textbf{Bim}}(\K,\K)$ with the category
$\R$.

Given small $K$-categories $\A$, $\B$ and $\C$, one has a functor
$$ -\otimes - :  \hskip .1in \text{\textbf{Bim}}(\B,\C) \times
\text{\textbf{Bim}}(\A,\B) \to \text{\textbf{Bim}}(\A,\C)$$ defined
as follows: Let $ {\phi}:\A \rightsquigarrow \B $ be an
$\A$-$\B$-bimodule and let $ {\varphi}: \B \rightsquigarrow \C $ be
a $\B$-$\C$-bimodule. For any given object $a \in \A$, the
$K$-functor $\phi(-, a): \B^{op} \to \R$ can be seen as a
$\K$-$\B$-bimodule ${\phi}(-, a): \K \rightsquigarrow \B $, while
for any $c \in \C$, the $K$-functor $\varphi(c, -): \B \to \R$ can
be seen as a $\B$-$\K$-bimodule ${\varphi}(c, -): \B\rightsquigarrow
\K $. Then the $\A$-$\C$-bimodule $\varphi \otimes \phi$ is the
$K$-functor defined by
$$(\varphi \otimes \phi) (c, a)=\varphi (c, -) \otimes
\phi (-, a), c \in \C, a \in \A.$$ This tensor product of bimodules
is associative up to canonical isomorphisms and admits as units the
$\B$-$\B$-bimodules $ 1_{\B}: \B \rightsquigarrow
 \B $ given by $1_{\B}(b, b')= \B(b, b')$; so that $1_{\B}
\otimes \phi \simeq \phi$ and $\varphi \otimes 1_{\B} \simeq
\varphi.$ Consequently, we have a bicategory $ \text{\textbf{Bim}}$
whose 0-cells are small $K$-categories, whose 1-cells are bimodules,
and whose 2-cells are $K$-functors between bimodules. The horizontal
product in $ \text{\textbf{Bim}}$ is the tensor product of
bimodules.

\bigskip
We shall need the following
\begin{lemma}\label{purepure}
\begin{itemize}
  \item [\emph{(i)}] A 2-cell $\tau: f \to f'$ in $\mathbf{Bim}$
  is right pure iff it is right $\K$-pure.
  \item [\emph{(ii)}] A 1-cell $f$ in $\mathbf{Bim}$ is right flat
  iff it is right $\K$-flat.

\end{itemize}
\end{lemma}
\begin{proof}(i). One direction is clear, so suppose that $\tau$ is right
$\K$-pure, and consider an arbitrary 1-cell $f: \C \to \A$. For any
pair $(b,c)\in \B \times \C,$ the $(b,c)$-component $(\tau \otimes
f)_{(b,c)}$ of the 2-cell $\tau \otimes f: f \otimes f \to f'
\otimes f $ is $\tau_{b,-} \otimes f(-,c)$, and identifying $f(-,c)$
with a bimodule $\K \to \A$, $\tau_{b,-} \otimes f(-,c)$ can be seen
as the $b$-component $(\tau \otimes f(-, c))_b$ of the 2-cell $\tau
\otimes f(-, c):f \otimes f(-, c)\to f' \otimes f(-, c)$. But since
$f$ is right $\K$-pure, $\tau \otimes f(-, c)$ is a monomorphism for
all $c \in \C$. Now, since any natural transformation between
$K$-functors is a monomorphism in the corresponding functor category
iff it is a componentwise monomorphism, it follows that $(\tau
\otimes f(-, c))_b$, and hence $\tau_{b,-} \otimes f(-,c)$, is a
monomorphism for all $b\in \B$ and $c \in \C$. Thus $\tau \otimes f$
is a monomorphism, proving that $f$ is right pure.

In a similar way one can show that a 1-cell $f$ in
$\text{\textbf{Bim}}$ is right flat iff it is right $\K$-flat.
\end{proof}

There is of course a dual result:
\begin{lemma}\label{purepuredual}
\begin{itemize}
  \item [\emph{(i)}] A 1-cell $f$ in $\mathbf{Bim}$ is left flat iff
  it is left $\K$-flat.
  \item [\emph{(ii)}] A 2-cell $\tau: f \to f'$ in $\mathbf{Bim}$
  is left pure iff it is left $\K$-pure.

\end{itemize}
\end{lemma}

Since the category $\R$, and hence any functor category $[\A,\R]$,
is equipped with a representative choice of (regular) monomorphisms,
we usually write $I$ rather than $[(I, i_I)]$ for any subobject of
an arbitrary functor $\A \to \R$.

Let now $\A$ be a small $K$-category and $\mathbb{S}=(S, e_S, m_S)$
be an $\A$-ring. An $\mathcal{A}$--subbimodule of $\mathbb{S}$ is a
subobject $I$ of $S(-, -): \A^{\text{op}}\to \R$ in the category
$[\A^{\text{op}}\otimes \A, \R] $. It can be regarded as the union
of all the sets $I(a,a')$ of morphisms of $\A$ such that, if $f,g: a
\to a'$ are in $I(a,a')$ and if $x: b \to a$ and $y: a' \to b'$ are
arbitrary morphisms in $\A$, then $y(f+g)x$ is in $I(b,b')$. Quite
obviously, $S$ is an $\mathcal{A}$--subbimodule of $\mathbb{S}$.
Moreover, it is easy to see that if $I, \, J$ are
$\mathcal{A}$--subbimodules of $\mathbb{S}$, then their product $IJ$
defined by
$$IJ(a,a')=\{\sum_{k \in K} f_k g_k, f_k \in I(a_k, a'), g_k \in J(a,a_k ),K
\,\,\text{is \,\,a\,\,finite\,\,set\,\,and\,\,} a, a_k \in \A \}$$
is also an $\mathcal{A}$--subbimodule of $\mathbb{S}$. With respect
to this product, the $\mathcal{A}$--subbimodules of $\mathbb{S}$
form a monoid and  $S$ is a two sided unit for this product. We
write $\textbf{I}_\A(\mathbb{S})$ for this monoid.

\begin{remark}\label{IJ}
One can check easily that $IJ$ is just the image of the composite
$$\xymatrix{I\otimes J \ar[r]^{i_I \otimes i_J}& S \otimes S \ar[r]^-{m_S}& S\,,}$$
where $i_I: I \to S$ and $i_J: J \to S$ are the canonical
embeddings.
\end{remark}

Let $\textbf{I}^l_\A(\mathbb{S})$ (resp.
$\textbf{I}^r_\A(\mathbb{S})$) denote the submonoid of
$\textbf{I}_\A(\mathbb{S})$ consisting of those
$\mathcal{A}$--subbimodules $I$ for which the multiplication
$$S \w I \to S, \,\, s \w i \to si$$
$$(\text{resp.}\,\, I \w S \to s, \,\, i \w s \to is) $$ is an isomorphism.

\begin{proposition}\label{iso0}
Let $\eta_f,\e_f : f \dashv f^* : \B \to \A$ be an adjunction in
$\text{\textbf{\emph{Bim}}}$. If the functor $f \otimes - :
[\A^{\text{op}}, \R] \to [\B^{\text{op}}, \R]$ is precomonadic, then
for all $(I,i_I), (J, i_J) \in \mathbf{I}^{\,l}_{\mathbf{Bim}(\A,
\A)}({S}_f)$, the canonical morphism $I \otimes J \to IJ$ is an
isomorphism.
\end{proposition}
\begin{proof} Consider the composite $$\xymatrix{I\otimes J \ar[r]^-{i_I \otimes J}&
S_f \otimes J \ar[r]^-{S_f \otimes i_J}&S_f \otimes S_f
\ar[r]^-{m_f}& S_f\,.}$$ Since the functor $$f \otimes - :
[\A^{\text{op}}, \R] \to [\B^{\text{op}}, \R]$$ (which can be seen
as the functor
$$\text{\textbf{Bim}}(\A, f)=f \w - : \text{\textbf{Bim}}(\K, \A)
\to \text{\textbf{Bim}}(\K, \B))$$ is assumed to be precomonadic, it
follows  that $\eta_f : 1 \to S_f$ is right $\K$-pure, and hence
right pure by Lemma \ref{purepure}. Thus, the unit $\eta_f : 1 \to
S_f$ of the $\text{\textbf{Bim}}(\A, \,\A)$-monoid $S_f$ is right
pure in the monoidal category $\text{\textbf{Bim}}(\A, \,\A)$. It
then follows from Proposition \ref{purepures} that $i_I \otimes J$
is a monomorphism. Moreover, since $(J, i_J) \in
\textbf{I}^{\,l}_{\text{\textbf{Bim}}(\A, \, \A)}(S_f)$,
\,$\xi^l_J=m_f \cdot (S_f \otimes i_J)$ is an isomorphism.
Consequently, the composite $m_f \cdot (i_I \otimes i_J)$ is a
monomorphism, and thus its image -which is just $IJ$ by Remark
\ref{IJ} - is isomorphic to $I \otimes J$.

\end{proof}
\begin{corollary}\label{iso1} In the situation of Proposition \ref{iso0}, the assignment
$$I \longrightarrow (I, i_I)$$ yields an isomorphism
$$\textbf{\emph{I}}^{\,l}_\A({S}_f) \simeq
\textbf{\emph{I}}^{\,l}_{\mathbf{Bim}(\A, \A)}({S}_f)$$ of monoids.
\end{corollary}

Dually, one has
\begin{corollary}\label{iso2} In the situation of Proposition \ref{iso0}, the assignment
$$I \longrightarrow (I, i_I)$$ yields an isomorphism
$$\textbf{\emph{I}}^{\,r}_\A({S}_f) \simeq
\textbf{\emph{I}}^{\,r}_{\mathbf{Bim}(\A, \A)}({S}_f)$$ of monoids.
\end{corollary}

\bigskip
Considering the composites
$$\xymatrix{\mathbf{I}^{\,l}_\A({S}_f) \simeq
\mathbf{I}^{\,l}_{\text{\textbf{Bim}}(\A, \A)}({S}_f)
\ar[r]^-{\Gamma_f}& \text{End}_\B(\mathfrak{C}_f,\mathfrak{C}_f)}$$
and
$$\xymatrix{(\textbf{I}^{\,r}_\A({S}_f))^{\text{op}} \simeq
(\mathbf{I}^{\,r}_{\text{\textbf{Bim}}(\A, \A)}({S}_f))^{\text{op}}
\ar[r]^-{\Gamma_{f^*}}&
\text{End}_\B(\mathfrak{C}_f,\mathfrak{C}_f)}$$$-$ which we still
call $\Gamma_f$ and $\Gamma_{f^*}$, respectively$-$ and combining
Theorems \ref{maincomonadic}, \ref{maincomonadicdual} and
\ref{separableiso} and Corollaries \ref{iso1} and \ref{iso2}, we get

\begin{theorem}\label{comonadicbim}
Let If $\eta_f,\e_f : f \dashv f^* : \B \to \A$ be an adjunction in
$\text{\textbf{\emph{Bim}}}$.
\begin{itemize}
\item [\emph{(i)}] If the functor $f \otimes - : [\A^{\text{op}},
\R] \to [\B^{\text{op}}, \R]$ is comonadic, the map $$\Gamma_f:
\textbf{\emph{I}}^{\,l}_\A({S}_f) \to
\text{\emph{End}}_\B(\mathfrak{C}_f,\mathfrak{C}_f)$$ is an
isomorphism of monoids.
  \item [\emph{(ii)}] If the functor $- \otimes f^* : [\A, \R] \to [\B, \R]$
  is comonadic, the map $$\Gamma_{f^*}:(\textbf{\emph{I}}^{\,r}_\A({S}_f))^{\text{op}}
  \to \text{\emph{End}}_\B(\mathfrak{C}_f,\mathfrak{C}_f) $$ is an isomorphism of monoids.
\end{itemize}
\end{theorem}

Now, combining Theorem \ref{comonadicbim}  with  Theorems
\ref{flatpureiso} and \ref{separableiso} and Lemmas \ref{purepure}
and \ref{purepuredual} we get:

\begin{theorem}\label{flatisobim}
Let $\eta_f,\e_f : f \dashv f^* : \B \to \A$ be an adjunction in
$\text{\textbf{\emph{Bim}}}$.
\begin{itemize}
\item [\emph{(i)}]If $f$ is right flat and $\eta_f$ is right pure,
then the map $$\Gamma_f: \textbf{\emph{I}}^{\,l}_\A({S}_f) \to
\text{\emph{End}}_\B(\mathfrak{C}_f,\mathfrak{C}_f)$$ is an
isomorphism of monoids.

\item [\emph{(ii)}]If $f$ is left flat and $\eta_f$ is left pure,
then the map
$$\Gamma_{f^*}:(\textbf{\emph{I}}^{\,r}_\A({S}_f))^{\text{op}}
  \to \text{\emph{End}}_\B(\mathfrak{C}_f,\mathfrak{C}_f) $$ is an isomorphism of monoids.

\item [\emph{(iii)}] If $f$ is a separable bimodule, then each of
the maps $\Gamma_f$ and $\Gamma_{f^*}$, is an isomorphism of
monoids.
\end{itemize}
\end{theorem}

Recall that a $K$-category $\A$ with only one object amounts to a
$K$-algebra $A$ (with unit) and that the category $\mathbf{Bim}(\K,
\, \A)$ is (isomorphic to) the category of right $A$-modules
$\text{Mod}_A$, while the category $\mathbf{Bim}(\A, \, \K)$ is
(isomorphic to) the category of left $A$-modules $_A\text{Mod}$.
Recall also that if $\A$ and $\B$ are two such $K$-categories with
corresponding $K$-algebras $A$ and $B$, respectively, then to give
an $\A$-$\B$-bimodule $ \MM: \A \rightsquigarrow \B $ is to give a
$A$-$B$-bimodule $M$ and that the diagrams
$$\xymatrix{\mathbf{Bim}(\K, \, \A) \ar@{=}[d]\ar[r]^-{\MM \otimes
-}& \mathbf{Bim}(\K, \, \B) \ar@{=}[d]\\
\text{Mod}_A \ar[r]_{- \otimes_A M}& \text{Mod}_B}$$ and
$$\xymatrix{\mathbf{Bim}(\B, \, \K)
\ar@{=}[d]\ar[r]^-{- \otimes \MM}& \mathbf{Bim}(\A, \, \K) \ar@{=}[d]\\
_B\text{Mod} \ar[r]_{M \otimes_B -}& _A\text{Mod}\,,}$$ where the
vertical morphisms are the isomorphisms, are both commutative.
Recall finally that, according e.g. to \cite{BD}, $\MM$ has a right
adjoint $\MM^*$ in $\mathbf{Bim}$ iff $M_B$ is finitely generated
and projective. When this is the case, the $\B$-$\A$-bimodule
$\MM^*: \B \to \A$ corresponds to the $B$-$A$-bimodule
$M^*=\text{Mod}_B (M, B)$. Moreover, the diagram
$$\xymatrix{\mathbf{Bim}(\K, \, \B)
\ar@{=}[d]\ar[r]^-{\MM^* \otimes
-}& \mathbf{Bim}(\K, \, \A) \ar@{=}[d]\\
\text{Mod}_B \ar[r]_{- \otimes_B M^*}& \text{Mod}_A}$$ commutes. It
follows that the $\A$-ring $S_{\MM}=\MM^* \otimes \MM$ corresponds
to the $A$-ring $S_M=M \otimes_B M^* \simeq \text{Mod}_B (M,M)$,
while the $\B$-coring $\mathfrak{C}_{\MM}=\MM \otimes \MM^*$
corresponds to the so called \emph{comatrix coring} $M^* \otimes_A
M$ corresponding to the $A$-$B$-bimodule $M$ (see \cite{EG1}). This
means that if $ \{(e_i, \, e^*_i)_{1\leq i \leq n}\}\subset M \times
M^*$ is a fixed dual basis for $M_B$, then $M^* \otimes_A M$ is a
$B$-coring with the following comultiplication $\Delta$ and and
counit $\varepsilon$
$$\Delta (m \otimes_A m^*)=\sum_i m^* \otimes_A e_i \otimes_B e^*_i \otimes _A m_i, \,\,
\varepsilon (m \otimes_A m^*)=m^*(m).$$

Since the horizontal composition in $\mathbf{Bim}$ is the opposite
of the usual tensor product of bimodules, we get that
$\mathbf{I}^{\,l}_A(S_{\mathcal{M}}) =
\mathbf{I}^{\,r}_A(S_M)^{op}$, and
$\mathbf{I}^{\,r}_A(S_{\mathcal{M}}) = \mathbf{I}^{\,
l}_A(S_M)^{op}$, where $\mathbf{I}^{\,l}_A(S_M)$ and $\mathbf{I}^{\,
r}_A(S_M)$ are the monoids associated to the ring homomorphism $A
\rightarrow S_M$ according to \cite{Masuoka:1989}. Applying now
Theorem \ref{comonadicbim} gives

\begin{theorem}\emph{(\cite{Me1})} Let $A, \, B$ be $K$-algebras, $M$ be an
$A$-$B$-bimodule  with $M_B$ finitely generated projective. Then:

\begin{itemize}

\item[\emph{(i)}] If the functor $- \otimes_A \!M :
\emph{\text{Mod}}_A \to \emph{\text{Mod}}_B$ is comonadic, then the
map
$$\Gamma_{M}:(\textbf{\emph{I}}^{\,r}_A(S_M))^{\text{op}}
  \to \emph{\text{End}}_{B-\text{cor}}(\mathfrak{C}_M)$$ that takes
$(I, i_I)\in \textbf{\emph{I}}^{\,r}_A(S_M)$ to the composite
$$\xymatrix{M^*\! \otimes_A \!M \ar[rr]^-{M^* \!\otimes_A \xi^{-1}_{i_I}}&&
M^*\! \otimes_A \!I \!\otimes_A \!M \ar[rr]^-{M^* \otimes_Ai_I
\otimes_A M} && M^*\! \otimes_A \! M \! \otimes_B \!M^* \otimes_A
\!M \ar[rr]^-{M^* \!\otimes_A \! M^* \! \otimes_B  \varepsilon}&&
M^* \! \otimes_A \!M}$$ is an isomorphism of monoids.

\item[\emph{(ii)}] If the functor $M^*\! \otimes_A - :\, _A
\emph{\text{Mod}} \to _B \!\emph{\text{Mod}}$ is comonadic, then the
map
$$\Gamma_{M^*}:\textbf{\emph{I}}^{\,l}_A(S_M)
  \to \emph{\text{End}}_{B-\text{cor}}(\mathfrak{C}_M)$$ that takes
$(I, i_I)\in \textbf{\emph{I}}^{\,r}_A(S_M)$ to the composite
$$\xymatrix{M^*\! \otimes_A \!M \ar[rr]^-{\xi^{-1}_{i_I} \!\otimes_A M}&&
M^*\! \otimes_A \!I \!\otimes_A \!M \ar[rr]^-{M^* \otimes_Ai_I
\otimes_A M} && M^*\! \otimes_A \! M \! \otimes_B \!M^* \otimes_A
\!M \ar[rr]^-{M^* \!\otimes_A \! M^* \! \otimes_B  \varepsilon}&&
M^* \! \otimes_A \!M}$$ is an isomorphism of monoids.

  \item[\emph{(iii)}] If $_A M_B$ is a separable bimodule, then
the maps $\Gamma_M$ and $\Gamma_{M^*}$ are both isomorphisms of
monoids.
\end{itemize}
\end{theorem}

From Theorem \ref{flatisobim} one obtains:

\begin{theorem}\cite[Theorem 2.5]{EG} Let $A, \, B$ be $K$-algebras, $M$ be an
$A$-$B$-bimodule  with $M_B$ finitely generated projective.

\begin{itemize}

\item[\emph{(i)}] If ${_A M}$ is faithfully flat, then the map
$\Gamma_{M}$ is an isomorphism of monoids.

\item[\emph{(ii)}] If ${M^*_A}$ is faithfully flat, then the map
$\Gamma_{M^*}$ is an isomorphism of monoids.

\item[\emph{(iii)}] If ${_A M_B}$ is a separable bimodule, then
the maps $\Gamma_M$ and $\Gamma_{M^*}$ are both isomorphisms of
monoids.
\end{itemize}

\end{theorem}

Let now $A$ be a ring, which is not assumed to have a unit and let
$\text{MOD}_A$ denote the category of all (not necessarily unital)
right $A$-modules. Following \cite{Quillen:notes}, a right
$A$-module $M$ is said to be \emph{firm} if the map
$$M \otimes_A A \to M , \,\, m \otimes_A a \to ma$$ is an
isomorphism. We write $\text{Mod}_A$ for the full subcategory of
$\text{MOD}_A$  whose objects are all right firm $A$-modules. In the
same way one defines the categories $_A\text{MOD}$ and $_A
\text{Mod}$. An $A$-$B$-bimodule is firm if it is firm both as a
left $A$-module and right $B$-module. A ring $A$ is firm if it is a
left (equivalently right) firm $A$-module. It is well-known that
when $A$ is a firm ring, the canonical embedding $i: \text{Mod}_A
\to \text{MOD}_A$ has a right adjoint functor $j : \text{MOD}_A \to
\text{Mod}_A$ that takes a (right) $A$-module $M$ to $M \otimes_A
A$. An object $M \in \text{MOD}_A$ is said \emph{closed} if the
canonical homomorphism $M \to \text{MOD}_A(A, M)$ given by $m
\longrightarrow (r \to mr)$ is an isomorphism. For any $M \in
\text{MOD}_A,$ the right $A$-module $\text{Mod}_A(A, M)$ is closed
and the assignment $M \to \text{Mod}_A(A, M)$ yields a functor $k:
\text{Mod}_A \to \text{MOD}_A$ which is right adjoint to the
reflection functor $j$ and which identifies the category of right
firm $A$-modules with the category of closed right $A$-modules.
Since the full subcategory of closed right $A$-modules is a
localization of $\text{MOD}_A$, it follows that $\text{Mod}_A$,
being equivalent to the category of right closed $A$-modules, is
also abelian. Note that kernels in $\text{Mod}_A$ cannot in general
be computed in $\text{MOD}_A$.

We write $\text{\textbf{Firm}}$ for the bicategory whose objects are
firm rings, whose 1-arrows are firm bimodules and whose 2-arrows are
homomorphism of firm bimodules (see \cite{GV}). The horizontal
composition is here the opposite of the usual tensor product of
bimodules.

Recall that a bicategory $\mathbb{B}$ is \emph{biclosed} if for any
1-cell $f: \A \to \B$ and any 0-cell $\C$, $$f \otimes - :
\mathbb{B}(\C, \A)\to \mathbb{B}(\C, \B)$$ has a right adjoint
$$[f, -]:\mathbb{B}(\C, \B)\to \mathbb{B}(\C, \A)$$ and $$- \otimes f : \mathbb{B}(\B, \C)
\to \mathbb{B}(\A, \C)$$ has a right adjoint $$\{f,
-\}:\mathbb{B}(\A, \C) \to \mathbb{B}(\B, \C). $$

\begin{proposition}The bicategory $\mathbf{Firm}$ is
biclosed.
\end{proposition}
\begin{proof} We have to show that for all firm rings $A$, $B$ and $C$ and
an arbitrary firm $A$-$B$-bimodule $M$, both functors
$$-\otimes_A \!M: \, _C \text{Mod}_A =\text{\textbf{Firm}}(C,A) \to
\text{\textbf{Firm}}(C,B)=_C \!\text{Mod}_B$$ and
$$M \otimes_B -: \, _B \text{Mod}_C =\text{\textbf{Firm}}(B,C) \to
\text{\textbf{Firm}}(A,C)=_A \!\text{Mod}_C$$ have a right adjoint.

According to Proposition 2.6 in \cite{GV}, for arbitrary two firm
rings $X$ and $Y$, the canonical isomorphism of categories
$t_{X,\,Y}:\, _X \!\text{MOD}_Y \to \text{MOD}_{X^{\text{op}}
\otimes Y} $ restricts to an isomorphism $\overline{t}_{X,\,Y}:\, _X
\!\text{Mod}_Y \to \text{Mod}_{X^{\text{op}} \otimes Y}$. Thus the
diagram $$\xymatrix{_X \!\text{Mod}_Y \ar[dd]_{i}
\ar[rr]^-{\overline{t}_{X,\,Y}} && \text{Mod}_{X^{\text{op}} \otimes
Y} \ar[dd]^{i}\\\\
_X \!\text{MOD}_Y  \ar[rr]_-{t_{X,\,Y}}&& \text{MOD}_{X^{\text{op}}
\otimes Y},\,}$$ where the vertical arrows are the canonical
embeddings, commutes. Since the functor $$i:
\text{Mod}_{X^{\text{op}} \otimes Y} \to \text{MOD}_{X^{\text{op}}
\otimes Y}$$ has a right adjoint $$j: \text{MOD}_{X^{\text{op}}
\otimes Y} \to \text{Mod}_{X^{\text{op}} \otimes Y}$$ that takes $M
\in \text{MOD}_{X^{\text{op}} \otimes Y} $ to $M
\otimes_{X^{\text{op}} \otimes Y} (X^{\text{op}} \otimes Y )$, the
functor $$i:\, _X \!\text{Mod}_Y \to\, _X \!\text{MOD}_Y$$ also has
a right adjoint $$j:\, _X \!\text{MOD}_Y \to\, _X \!\text{Mod}_Y$$
that takes $M \in _X \!\text{MOD}_Y $ to $X \!\otimes \!M \otimes
\!Y.$

Let now $M \in \text{\textbf{Firm}}(A,B)$, and let $C$ be an
arbitrary firm ring. Since the functor $$-\otimes_A \!M: \, _C
\text{MOD}_A  \to _C \!\text{MOD}_B$$ admits as a right adjoint the
functor $$\text{MOD}_B (M, -) :\, _C \text{MOD}_B \to _C
\text{MOD}_A, $$ it follows from the commutativity of the following
diagram $$\xymatrix{_C \!\text{Mod}_A \ar[rr]^{- \otimes_A M}
\ar[dd]_{i}&& _C
\!\text{Mod}_B \ar[dd]^{i}\\\\
_C \!\text{MOD}_A \ar[rr]_{-\otimes_A M}&& _C \!\text{MOD}_B} $$
that the functor
$$-\otimes_A \!M: \, _C \text{Mod}_A \to _C \!\text{Mod}_B$$ admits
as a right adjoint the functor $$j \cdot \text{MOD}_B (M, -) \cdot i
= C \!\otimes_C\! \text{Mod}_B (M, -)\!\otimes_A \!A.$$

Symmetrically, one shows that the functor $$B \!\otimes_B\!
\text{Mod} (M, -)\!\otimes_C \!C:\, _A \text{Mod}_C \to _B
\!\text{Mod}_C$$ is right adjoint to $$M \otimes_B -: \, _B
\text{Mod}_C \to _A \!\text{Mod}_C.$$ This completes the proof.
\end{proof}

We shall need the following well-known characterization of left
adjoint 1-cells in a biclosed bicategory (see, for example,
\cite{K}):

\begin{theorem}For any 1-cell $f : \A \rightsquigarrow \B$ in a biclosed bicategory
$\mathbb{B}$, the following are equivalent:
\begin{itemize}

\item[\emph{(i)}] $f$ has a right adjoint $f^*$;

\item[\emph{(ii)}] the morphism $[f, 1_\B]\otimes f \to [f,f]$
that corresponds under the adjunction $$\eta, \varepsilon : f
\otimes - \dashv [f, -]: \mathbb{B}(\B, \B) \to \mathbb{B}(\B, \A)$$
to the composite $\xymatrix{f \otimes [f, 1_\B] \otimes f
\ar[r]^-{\varepsilon_f \otimes f}& 1_\B  \otimes f \ar[r]^-{\simeq}
& f}$, is an isomorphism.
\end{itemize}When these conditions hold, one has $[f, g]\simeq f^* \otimes
g$ for all 1-cells $g: \C \rightsquigarrow \B $. In particular, $f^*
\simeq [f, 1_\B]$.
\end{theorem}

From this general theorem, we deduce:

\begin{theorem}\label{adjointfirm}
For any 1-cell $M : \A \rightsquigarrow \B$ in $\mathbf{Firm}$, the
following are equivalent:
\begin{itemize}

\item[\emph{(i)}]$M^{\dag}=B \!\otimes_B\! \emph{\text{Mod}}_B (M,
B)\!\otimes_A \!A$ is right adjoint to $M$ in $\mathbf{Firm}$;

\item[\emph{(ii)}]for any $X \in \mathbf{Firm}(C, B),$ the map
$$X \!\otimes_B\! \mathrm{Mod}_B (M, B)\!\otimes_A \!A \to
\mathrm{Mod}_B (M, X)\!\otimes_A \!A$$$$(x \otimes_B f \otimes_A
a)\longrightarrow ((m \to xf(m))\otimes_A a) $$ is an isomorphism.

\item[\emph{(iii)}]the functor $\mathrm{Mod}_B (M, -)\!\otimes_A
\!A: \mathrm{Mod}_B \to \mathbf{Ab}$ preserves all small colimits
\end{itemize}
\end{theorem}
\begin{proof} $\text{(i)}$ and $\text{(ii)}$ are equivalent by the previous
theorem, while $\text{(ii)}$ implies $\text{(iii)}$ trivially. Now,
if the functor $\text{Mod}_B (M, -)\!\otimes_A \!A$ preserves all
small colimits, then since for any $X \in \text{\textbf{Firm}}(C,
B),$ $X \otimes_BB \simeq X$ and since $X \otimes_BB$ can be seen as
a small colimt, we have: $$\text{Mod}_B (M, X)\!\otimes_A \!A \simeq
\text{Mod}_B (M, X \otimes_BB)\!\otimes_A \!A \simeq X \!\otimes_B\!
\text{Mod}_B (M, B)\!\otimes_A \!A.$$ Thus $\text{(iii)}$ implies
$\text{(ii)}$.
\end{proof}

\begin{remark}If $_AM_B$ is such that $M_B$ is finitely
generated and projective, then it is not hard to see that $M$ has a
right adjoint in $\text{\textbf{Firm}}$.
\end{remark}

Let $A$ and $B$ be firm rings  and let $M$ be a firm
$A$-$B$-bimodule. If $M$ admits a right adjoint in
$\text{\textbf{Firm}}$, then we know from Theorem \ref{adjointfirm}
that its right adjoint $M^\dag$ is $B \!\otimes_B\! \text{Mod}_B (M,
B)\!\otimes_A \!A$. Then the corresponding $A$-ring is $$S_M= M
\!\otimes_B\! M^\dag= M\! \otimes_B B \!\otimes_B\! \text{Mod}_B (M,
B)\!\otimes_A \!A \simeq M \!\otimes_B\!\text{Mod}_B (M,
B)\!\otimes_A \!A,$$ while the corresponding $B$-coring is
$$\mathfrak{C}_M= M^\dag\! \otimes_A \!M =B \!\otimes_B\! \text{Mod}_B (M,
B)\!\otimes_A \!A \! \otimes_A \!M \simeq B \!\otimes_B\!
\text{Mod}_B (M, B)\!\otimes_A \!M.$$

This coring was considered in \cite{GT:comonad}, and it generalizes
the (infinite) comatrix corings considered in \cite{EGT:infinite},
for $A$ with enough idempotents, and $B$ unital, and in \cite{GTV},
for $A$ firm and $B$ unital.

Let $\mathcal{M}$ (resp. $\mathcal{M}^\dag$) denote for a while the
firm $A$--$B$--bimodule $M$  (resp. the firm $B$--$A$--bimodule
$M^\dag$) when considered in the bicategory $\mathbf{Firm}$. Arguing
as in the proof of Proposition \ref{iso0}, one sees that if the
functor
$$\mathcal{M} \otimes - = -\otimes_A M : \text{Mod}_A \to \text{Mod}_B$$ (resp. the
functor $$- \otimes \mathcal{M}^\dag = M^\dag\! \otimes_A -\simeq B
\!\otimes_B\! \text{Mod}_B (M, B)\!\otimes_A -: \, _A\!\text{Mod}
\to _B\! \text{Mod})$$ is precomonadic, then the assignment $I \to
(I, i_I)$ yields an isomorphism of monoids
$$\mathbf{I}^r_A(S_M)^{op} \simeq \mathbf{I}^l_{\text{\textbf{Firm}}(A,\,A)}(S_{\mathcal{M}})$$
$$(\text{resp.} \,\,\mathbf{I}^l_A(S_M)^{op} \simeq \mathbf{I}^r_{\text{\textbf{Firm}}(A,\,A)}(S_\mathcal{M})),$$
where $\mathbf{I}^{r}_A(S_M)$ (resp. $\mathbf{I}^{S}_A(S_M)$)
consists of the firm $A$--subbimodules of $S_M$ such that the
multiplication map $I \otimes_A S_M \rightarrow S_M$ (resp. $S_M
\otimes_A I \rightarrow S_m$) is an isomorphism, which turns out to
be a monoid with the product given by the usual multiplication of
subsets of $S_M$.
\medskip

 Applying now Theorems \ref{maincomonadic} and
\ref{maincomonadicdual}, and Proposition \ref{separable}, gives:

\begin{theorem} Let $A, \, B$ be firm rings, $M$ be a firm
$A$-$B$-bimodule  that admits a right adjoint in $\mathbf{Firm}$.
Then

\begin{itemize}

\item[\emph{(i)}] if the functor $- \otimes_A \!M :
\emph{\text{Mod}}_A \to \emph{\text{Mod}}_B$ is comonadic, then the
map
$$\Gamma_{M}:(\textbf{\emph{I}}^{\,r}_A(S_M))^{\text{op}}
  \to \emph{\text{End}}_{B-\text{cor}}(\mathfrak{C}_M)$$ is an isomorphism of
  monoids;

\item[\emph{(ii)}] if the functor $B \!\otimes_B\! \emph{\text{Mod}}_B (M, B)\!\otimes_A -:\,
_A \emph{\text{Mod}} \to _B \!\emph{\text{Mod}}$ is comonadic, then
the map
$$\Gamma_{M^\dag}:\textbf{\emph{I}}^{\,l}_A(S_M)
  \to \emph{\text{End}}_{B-\text{cor}}(\mathfrak{C}_M)$$ is an isomorphism of monoids.

\item[\emph{(iii)}] If $_A M_B$ is a separable firm bimodule, then
the maps $\Gamma_M$ and $\Gamma_{M^\dag}$ are both isomorphisms of
monoids.
\end{itemize}
\end{theorem}
\bibliographystyle{amsplain}

\end{document}